\documentclass[12pt]{article}
\usepackage{enumerate}
\usepackage{amsfonts}
\usepackage{amsmath}
\usepackage{amssymb}

\newcounter{defin}  \newcounter{lemma}  \newcounter{theorem}
\newcounter{property} \newcounter{corol}  \newcounter{remark} \newcounter{example}

\newenvironment{lemma}{\par\refstepcounter{lemma}
     \textbf{Lemma \thelemma.} }{\rm\par}
\newenvironment{theorem}{\par\refstepcounter{theorem}
     \textbf{Theorem \thetheorem.}\ }{\rm\par}
\newenvironment{property}{\par\refstepcounter{property}
     \textbf{Proposition \theproperty.}\ }{\rm\par}
\newenvironment{corollary}{\par\refstepcounter{corol}
     \textbf{Corollary \thecorol.} }{\rm\par}
\newenvironment{definition}{\par\refstepcounter{defin}
     \textbf{Definition \thedefin.}\ }{\rm\par}
\newenvironment{remark}{\par\refstepcounter{remark}
     \textbf{Remark \theremark.}}{\rm\par}
\newenvironment{example}{\par\refstepcounter{example}
     \textbf{Example \theexample.}}{\rm\par}

\begin{document}
\title{Continuity condition for concave functions on convex $\mu$-compact sets and its applications in quantum physics}

\author{M.E.~Shirokov\thanks{e-mail:msh@mi.ras.ru}\\\\
Steklov Mathematical Institute, Moscow, Russia}

\date{} \maketitle\vspace{-15pt}

\begin{abstract}
A method of proving local continuity of concave functions on convex
set possessing the $\mu$\nobreakdash-\hspace{0pt}compactness property is presented. This
method is based on a special approximation of these functions.

The class of $\mu$\nobreakdash-\hspace{0pt}compact sets can be considered as a natural
extension of the class of compact metrizable subsets of locally
convex spaces, to which particular results well known for
compact sets can be generalized.

Applications of the obtained continuity conditions to analysis of
different entropic characteristics of quantum systems and channels are
considered.
\end{abstract}\vspace{-15pt}

\tableofcontents

\section{Introduction}

A problem of proving local continuity of a given concave (convex)
function defined on a given convex set arises naturally in different
fields of mathematics. For example, in mathematical physics this
problem appears in analysis of entropy-type functions on
a set of states of some physical system.
In some cases this problem can be solved by using  general results of convex analysis,
but sometimes it is difficult to apply them to a given function defined
on a convex set not satisfying particular requirements (compactness, existence of inner points, etc.)

In this paper we consider local continuity conditions for concave
functions on $\mu$\nobreakdash-\hspace{0pt}compact convex sets based on a special
approximation of these functions. The class of $\mu$\nobreakdash-\hspace{0pt}compact sets
(see Definition \ref{mu-comp} in Section 2) can be considered as a
natural extension of the class of compact metrizable subsets of
locally convex spaces, to which particular results well known
for compact sets can be generalized \cite{P&Sh}. This class contains
all compact sets as well as many noncompact sets widely used in
applications. The simplest examples of noncompact $\mu$\nobreakdash-\hspace{0pt}compact
convex sets are the positive part of the unit ball of the Banach
space $\ell_1$ and its closed subset consisting of all countable probability distributions. Other examples
and simple criteria of the $\mu$\nobreakdash-\hspace{0pt}compactness property can be found
in \cite{P&Sh}.

For applications in quantum physics it is essential that the convex
set of  positive operators in a separable Hilbert space with unit
trace, generally called quantum states, is $\mu$\nobreakdash-\hspace{0pt}compact. In fact,
it is necessity to explore continuity properties of several entropic
characteristics of quantum states, in particular, of the von Neumann
entropy, that provides a basic impetus to find universal method of proving
local continuity of these characteristics. In \cite{Sh-11} this method was developed by using some special properties of the set
of quantum states. In this paper we show that it can be generalized to the class of $\mu$\nobreakdash-\hspace{0pt}compact
convex sets by using slightly different argumentation.

The paper is organized as follows. In Section 2 notations and
basic results used in the subsequent sections are presented. In
Section 3 we consider several properties of a convex set following
from its $\mu$\nobreakdash-\hspace{0pt}compactness and stability (see Definition
\ref{stability} in Section 2). Section 4 is devoted to a special
approximation technic for concave functions. In Section 5 the
continuity conditions based on this technic are presented.
Applications to quantum physics extending the results of
\cite{Sh-11} are considered in Section 6.

\section{Basic notations}

In what follows $\mathcal{A}$~is a bounded convex complete separable
metrizable subset of some locally convex space.\footnote{This means
that the topology on the set $\mathcal{A}$ is defined by a countable
subset of the family of seminorms, generating the topology of the
entire locally convex space, and this set is separable and complete
in the metric generated by this subset of seminorms.} The set of
extreme points of the set $\mathcal{A}$ will be denoted
$\mathrm{extr}(\mathcal{A})$.

Let $\mathrm{cl}(\mathcal{B})$, $\mathrm{co}(\mathcal{B})$,
$\sigma\textrm{-}\mathrm{co}(\mathcal{B})$ and
$\overline{\mathrm{co}}(\mathcal{B})$ be respectively  the closure,
the convex hull, the $\sigma$\nobreakdash-\hspace{0pt}convex hull\footnote{$\sigma\textrm{-}\mathrm{co}(\mathcal{B})$
is the set of all countable convex combinations of points in
$\mathcal{B}$.}
and the convex closure of a subset $\mathcal{B}\subseteq\mathcal{A}$
\cite{J&T,R}.

For an arbitrary closed subset $\mathcal{B}\subseteq\mathcal{A}$
denote by $C(\mathcal{B})$ the set of all continuous bounded
functions on $\mathcal{B}$, denote by $M(\mathcal{B})$ and
$M^a(\mathcal{B})$ respectively the set of all Borel probability
measures on $\mathcal{B}$ and its subset consisting of atomic
measures. We always assume that the set $M(\mathcal{B})$ and arbitrary its subsets are endowed
with the weak convergence topology~\cite{Bil,Par}.

With an arbitrary measure $\mu\in M(\mathcal{B})$ we associate its
barycenter (average)~$\mathbf b(\mu)\in
\overline{\mathrm{co}}(\mathcal{B})$, which is defined by the Pettis
integral~(see \cite{Ali,V&T})
\begin{equation}\label{eq1}
\mathbf b(\mu)\ =\ \int_{\mathcal{B}}x \mu(dx).
\end{equation}

If $\mu$ is a measure in $M^a(\mathcal{B})$ "consisting" of atoms
$\{x_i\}$ with the corresponding weights $\{\pi_i\}$ then  $\mathbf
b(\mu)=\sum_i \pi_i x_i$. The above measure  will be denoted $\sum_i
\pi_i\delta(x_i)$ or, briefly, $\{\pi_i, x_i\}$.\smallskip

For a given Borel function $f$ on a closed subset $\mathcal{B}\subseteq\mathcal{A}$ consider the
functional
\begin{equation}\label{functional}
 M(\mathcal{B})\ni\mu\mapsto \mathbf{f}(\mu)=\int_{\mathcal{B}}f(x)\mu(dx).
\end{equation}
It is easy to show that this functional is lower semicontinuous
(correspondingly, upper semicontinuous) provided the function $f$ is
lower semicontinuous and lower bounded (correspondingly, upper
semicontinuous and upper bounded) on the set $\mathcal{B}$
\cite{Bil}. \medskip

For arbitrary $x\in\overline{\mathrm{co}}(\mathcal{B})\,$ let
$M_{x}(\mathcal{B})$ and $M^a_{x}(\mathcal{B})$ be respectively
convex subsets of $M(\mathcal{B})$ and of $M^a(\mathcal{B})$
consisting of such measures  $\mu$  that $\,\mathbf b(\mu) =
x$.\smallskip

The barycenter map
\begin{equation}\label{b-map}
M(\mathcal{A})\ \ni \ \mu\; \mapsto \; \mathbf b(\mu)\ \in
 \mathcal{A}
\end{equation}
is continuous (this can be shown easily by applying Prokhorov's
theorem~\cite[Ch.II, Th.6.7]{Par}). Hence the image of any compact
subset of $M(\mathcal{A})$ under this  map is a compact subset of
$\mathcal{A}$. The $\mu$\nobreakdash-\hspace{0pt}compact sets are defined by the converse
requirement \cite{P&Sh}.\medskip

\begin{definition}\label{mu-comp}
A set $\mathcal{A}$ is called $\,\mu$-{\it compact } if the preimage
of any compact subset of $\mathcal{A}$ under barycenter map
(\ref{b-map}) is a compact subset of $M(\mathcal{A})$.
\end{definition}
\medskip

Any compact set is $\mu$\nobreakdash-\hspace{0pt}compact, since compactness of $\mathcal{A}$
implies compactness of $M(\mathcal{A})$~\cite{Par}. Properties
of $\mu$\nobreakdash-\hspace{0pt}compact sets are studied in detail in~\cite{P&Sh}, where
$\mu$\nobreakdash-\hspace{0pt}compactness of several noncompact sets widely used in
applications has been proved (for example, of the set of all Borel
probability measures on an arbitrary complete separable metric space
endowed with the weak convergence topology and of  the set of
quantum states -- density operators in a separable Hilbert
space).\smallskip

The $\mu$\nobreakdash-\hspace{0pt}compactness property of a convex set is not purely
topological but reflects a special relation between the topology
and the convex structure of this set.\smallskip

An another relation between the topology and the convex structure of
a convex set is expressed by the notion of (convex) stability
\cite{Stefania}.\smallskip

\begin{definition}\label{stability}
A set $\mathcal{A}$ is called \emph{stable} if the map $\
\mathcal{A}\times\mathcal{A}\, \ni\, (x,y)\ \mapsto
  \, \displaystyle\frac{x\, +\, y}{2}\, \in\, \mathcal{A}\,$ is
open.
\end{definition}\smallskip

The notion of stability of a convex subset of a linear topological
space appeared at the end of 1970's as a result of study of
convex compact sets, which leaded in particular to proving equivalence of
the following properties a convex compact set $\mathcal{A}$:
\begin{enumerate}[(i)]
 \item \textit{the set $\mathcal{A}$ is stable;}
 \item \textit{the map $M(\mathcal{A})\ni\mu\ \mapsto \ \mathbf b(\mu)\in \mathcal{A}\,$
is open;}
 \item \textit{the map $M\left(\overline{\mathrm{extr}\,
\mathcal{A}}\right)\ni\mu\ \mapsto \
 \mathbf b(\mu)\in \mathcal{A}\,$
is open;}
 \item \textit{the convex hull \footnote{The convex hull of a function is the maximal convex function majorized by this function \cite{J&T}.} of an arbitrary continuous function on $\mathcal{A}$ is continuous;}
 \item \textit{the convex hull of an arbitrary concave continuous function on $\mathcal{A}$ is continuous.}
\end{enumerate}
Essential parts of the above assertion was obtained by Vesterstrom
\cite{Ves}, its complete version was proved by O'Brien \cite{Brien}.
This assertion (called the Vesterstrom-O'Brien theorem in what follows) does not hold for noncompact convex sets in general,
but it can be extended to convex $\mu$\nobreakdash-\hspace{0pt}compact sets
\cite[Theorem 1]{P&Sh}.

In $\mathbb{R}^{2}$ stability holds for an arbitrary convex compact
set, in $\mathbb{R}^{3}$ it is equivalent to closedness of the set of
extreme points of a convex compact
set while in $\mathbb{R}^{n}, n>3,$ it is stronger than
the last property \cite{Brien}. A full characterization of the
stability property in finite dimensions is obtained in
\cite{Stefania}. In infinite dimensions stability is proved for the unit
ball in some Banach spaces and for the positive part of the
unit ball in Banach lattices in which the unit ball is stable \cite{Grzaslewicz}.

The simplest example of a noncompact $\mu$\nobreakdash-\hspace{0pt}compact convex stable set
is the set $\mathfrak{P}_{+\infty}$ of all probability distributions
with countable number of outcomes (considered as a subset of the
Banach space $\ell_1$). This is a partial case of the more general example -- the convex set of all Borel probability measures on any complete separable metric space endowed with the weak convergence topology. The $\mu$\nobreakdash-\hspace{0pt}compactness and stability of this set are established respectively in \cite[Corollary 4]{P&Sh} and in \cite[Theorem 2.4]{Eifler}.\medskip

We will use the following two strengthened versions of the notion of a concave function.\smallskip

A semibounded (upper or lower bounded) function $f$ on a convex set
$\mathcal{A}$ is called
$\sigma$\nobreakdash-\hspace{0pt}\textit{concave} if the discrete
Jensen's inequality
$$
f(\mathbf{b}(\{\pi_{i},x_{i}\}))\geq\sum_{i}\pi_{i}f(x_{i})
$$
holds for an arbitrary measure $\{\pi_{i},x_{i}\}$ in
$M^a(\mathcal{A})$.

A semibounded universally measurable\footnote{This means that the
function $f$ is measurable with respect to any measure in
$M(\mathcal{A})$.} function $f$ on a convex set $\mathcal{A}$ is
called $\mu$\nobreakdash-\hspace{0pt}\textit{concave} if the
integral Jensen's inequality
$$
f(\mathbf{b}(\mu))\geq\int_{\mathcal{A}}f(x)\mu(dx)
$$
holds for an arbitrary measure $\mu$ in $M(\mathcal{A})$.

$\sigma$\nobreakdash-\hspace{0pt}\textit{convexity} and
$\mu$\nobreakdash-\hspace{0pt}\textit{convexity} of a function $f$
are naturally defined via the above notions applied to the function
$-f$.

Examples of semibounded functions, which are convex but not
$\sigma$\nobreakdash-\hspace{0pt}convex or
$\sigma$\nobreakdash-\hspace{0pt}convex but not
$\mu$\nobreakdash-\hspace{0pt}convex, are considered in \cite[Section
3]{Sh-9}.

The following lemma contains sufficient conditions for
$\sigma$\nobreakdash-\hspace{0pt}concavity and
$\mu$\nobreakdash-\hspace{0pt}concavity of a concave function, which
can be proved easily (see the Appendix in \cite{Sh-9}).
\medskip

\begin{lemma}\label{yensen}
\emph{Let $f$ be a concave function on a convex set $\mathcal{A}$.}
\begin{enumerate}[A)]
    \item \emph{If $f$ is lower bounded then $f$ is
$\sigma$\nobreakdash-\hspace{0pt}concave.}
    \item \emph{If $f$ is either lower
semicontinuous and lower bounded or upper semicontinuous then $f$ is
$\mu$\nobreakdash-\hspace{0pt}concave.}
\end{enumerate}
\end{lemma} \medskip

\begin{remark}\label{perturb}
Assuming that the metric
$\mathbf{d}(\cdot\,,\cdot)$ on the set $\mathcal{A}$ is defined as follows
$$
\mathbf{d}(x,y)=\sum_{k=1}^{+\infty}2^{-k}\frac{\|x-y\|_k}{1+\|x-y\|_k}, \quad
x,y\in\mathcal{A},
$$
where $\{\|\cdot\|_k\}_{k=1}^{+\infty}$ is the countable
family of seminorms generating the
topology on this set,
it is easy to obtain the following estimation
$$
\mathbf{d}(\alpha x+(1-\alpha)y, \alpha'x'+(1-\alpha')y')\leq
2\delta+C_{x,y}(\varepsilon)
$$
valid for any $x,y,x',y'$ in $\mathcal{A}$ and any
$\alpha,\alpha'$ in $[0,1]$ such that
$\mathbf{d}(x,x')<\delta$, $\,\mathbf{d}(y,y')<\delta$ and
$|\alpha-\alpha'|<\varepsilon$, where
$C_{x,y}(\varepsilon)=\sum_{k=1}^{+\infty}2^{-k}\frac{\varepsilon\|x-y\|_k}{1+\varepsilon\|x-y\|_k}$ is a function such that $\lim_{\varepsilon\rightarrow+0}C_{x,y}(\varepsilon)=0$.
\end{remark}\medskip

\textbf{Note:}  In what follows continuity of a function $f$ on a
subset $\mathcal{B}\subset\mathcal{A}$ means continuity of the
restriction $f|_\mathcal{B}$ of the function $f$ to the subset $\mathcal{B}$, which implies
finiteness of this restriction (in contrast to lower or upper
semicontinuity).

\section{Some implications of $\mu$-compactness and stability of a convex set}

In this section we consider auxiliary  results  used in the main
part of the paper. \smallskip

We begin with  several simple lemmas. \smallskip

\begin{lemma}\label{bar-dec}
\textit{Let $\,\mathcal{B}$ be a closed subset of a convex
$\mu$\nobreakdash-\hspace{0pt}compact set $\mathcal{A}$. Then for
arbitrary $x_0$ in $\,\overline{\mathrm{co}}(\mathcal{B})$ there
exists a measure $\mu_0$ in $M(\mathcal{B})$ such that
$\,x_0=\mathbf{b}(\mu_0)$}.
\end{lemma}\smallskip

\textbf{Proof.} Let $x_{0}\in\overline{\mathrm{co}}(\mathcal{B})$
and $\{x_{n}\}\subset\mathrm{co}(\mathcal{B})$ be a sequence
converging to $x_{0}$. For each $n\in\mathbb{N}$ there exists  a
measure  $\mu_{n}\in M(\mathcal{B})$ with finite support such that
$x_{n}=\textbf{b}(\mu_{n})$. By $\mu$\nobreakdash-\hspace{0pt}compactness of the set
$\mathcal{A}$ the sequence $\{\mu_{n}\}$ has a partial limit
$\mu_{0}\in M(\mathcal{B})$. Continuity of the map
$\mu\mapsto\textbf{b}(\mu)$ implies $\textbf{b}(\mu_{0})=x_{0}$.
$\square$\medskip

\begin{lemma}\label{density}
\textit{Let $\mathcal{A}$ be a convex
$\mu$\nobreakdash-\hspace{0pt}compact set such that the set
$\,\mathrm{extr}\mathcal{A}$ is closed and
$\mathcal{A}=\sigma\textup{-}\mathrm{co}(\mathrm{extr}\mathcal{A})$.
Then an arbitrary measure $\mu_{0}$ in $M(\mathrm{extr}\mathcal{A})$
can be approximated by a sequence $\{\mu_{n}\}$ of measures in
$M^{a}(\mathrm{extr}\mathcal{A})$ such that
$\,\mathbf{b}(\mu_{n})=\mathbf{b}(\mu_{0})$ for all $\,n$.}
\end{lemma}\smallskip

\textbf{Proof.} Consider the Choquet ordering on the set
$M(\mathcal{A})$. We say that $\mu\succ\nu$ if and only if
$$
\int_{\mathcal{A}}f(x)\mu(dx) \geq\int_{\mathcal{A}}f(x)\nu(dx)
$$
for any convex continuous bounded function $f$ on the set
$\mathcal{A}$ \cite{Edgar}.\medskip

For a given measure $\mu_{0}$ in $M(\mathrm{extr}\mathcal{A})$ it is easy to construct
a sequence $\{\mu_{n}\}$ of measures in $M(\mathcal{A})$ with
finite support converging to the measure $\mu_{0}$ such that
$\mathbf{b}(\mu_{n})=\mathbf{b}(\mu_{0})$ for all $n$. Decomposing
each atom  of the measure $\mu_{n}$ into convex combination of
extreme points we obtain the measure $\hat{\mu}_{n}$ in
$M^{\mathrm{a}}(\mathrm{extr}\mathcal{A})$ with the same barycenter.
It is easy to see that $\hat{\mu}_{n}\succ\mu_{n}$. By
$\mu$\nobreakdash-\hspace{0pt}compactness of the set $\mathcal{A}$ the sequence
$\{\hat{\mu}_{n}\}_{n>0}$ is relatively compact. This implies
existence of subsequence $\{\hat{\mu}_{n_{k}}\}$ converging to a
measure $\hat{\mu}_{0}$ in $M(\mathrm{extr}\mathcal{A})$. Since $\hat{\mu}_{n_{k}}\succ\mu_{n_{k}}$
for all $k$, the definition of the weak convergence implies
$\hat{\mu}_{0}\succ\mu_{0}$ and hence $\hat{\mu}_{0}=\mu_{0}$ by
maximality of the measure $\mu_{0}$ with respect to the Choquet
ordering (which follows from coincidence of this ordering with the
dilation ordering \cite{Edgar}). $\square$
\medskip

\begin{lemma}\label{g-mu-comp}
\textit{Let $\mathcal{A}$ be a convex
$\mu$\nobreakdash-\hspace{0pt}compact set and
$\,\{\{\pi^{n}_{i},x^{n}_{i}\}_{i=1}^{m}\}_{n}$ be a sequence of
measures in $M^a(\mathcal{A})$ having $\,m<+\infty$ atoms such that
the sequence $\{\sum_{i=1}^{m}\pi^{n}_{i}x^{n}_{i}\}_{n}$ of their
barycenters converges to a point  $\,x_{0}\in\mathcal{A}$. There
exists a subsequence
$\{\{\pi^{n_k}_{i},x^{n_k}_{i}\}_{i=1}^{m}\}_{k}$ converging to a
particular measure\footnote{We do not assert that $x^{0}_{i}\neq
x^{0}_{j}$ for all $i\neq j$.} $\{\pi^{0}_{i},x^{0}_{i}\}_{i=1}^{m}$
with the barycenter $\,x_{0}$ in the following sense
$$
\lim_{k\rightarrow+\infty}\pi^{n_k}_{i}=\pi^{0}_{i}\quad\textrm{and}\quad\pi^{0}_{i}>0\;\,\Rightarrow\,
\lim_{k\rightarrow+\infty}x^{n_k}_{i}=x^{0}_{i},\quad
i=\overline{1,m}.
$$}
\end{lemma}

\textbf{Proof.} It is sufficient to note that $\mu$\nobreakdash-\hspace{0pt}compactness of
the set $\mathcal{A}$ implies relative compactness of the sequence
$\{\{\pi^{n}_{i},x^{n}_{i}\}_{i=1}^{m}\}_{n}$ and that the set of
measures having $m$ atoms is a closed subset of $M(\mathcal{A})$.
$\square$
\medskip

Let $\mathfrak{P}_{n}$ be the set of all probability distributions
with $n\leq+\infty$ outcomes. \medskip
\begin{lemma}\label{a-k-closed} \textit{Let $\,\mathcal{A}_1$ be a closed subset of a convex $\mu$\nobreakdash-\hspace{0pt}compact set $\,\mathcal{A}$.}
\begin{enumerate}[A)]
\item  \textit{The set
\begin{equation}\label{a-k}
\mathcal{A}_k=\left\{\left.\sum_{i=1}^{k}\pi_i x_i\, \right|
\{\pi_i\}\in\mathfrak{P}_k,\, \{x_i\}\subset\mathcal{A}_1\right\}
\end{equation}
is closed for each $\,k\in\mathbb{N}$.}
\item \textit{Let $f$ be a concave nonnegative function on the set
$\,\mathcal{A}$, which takes a finite value at least at one point in
$\,\mathcal{A}_1$. If this function is upper continuous on the set
$\,\mathcal{A}_k$ defined by (\ref{a-k}) for each $k$ then it is
bounded on the set $\,\mathcal{A}_k$ for each $k$.}
\end{enumerate}
\end{lemma}\smallskip

\textbf{Proof.} A) This assertion directly follows from Lemma
\ref{g-mu-comp}.

B) Suppose there exists a sequence $\{x_n\}\subset\mathcal{A}_k$
such that $\lim_{n\rightarrow+\infty}f(x_n)=+\infty$. Let $y_0$ be
a point in $\mathcal{A}_1$ with finite $f(y_0)$. Consider the
sequence $\{\lambda_n
x_n+(1-\lambda_n)y_0\}\subset\mathcal{A}_{k+1}$, where
$\lambda_n=1/f(x_n)$. This sequence converges to the point $y_0$ (since the set $\mathcal{A}$ is bounded),
but concavity of the function $f$ implies
$$
\liminf_{n\rightarrow+\infty}f(\lambda_n x_n+(1-\lambda_n)y_0)\geq
\liminf_{n\rightarrow+\infty}\left(\lambda_nf(x_n)+(1-\lambda_n)
f(y_0)\right)=1+f(y_0),
$$
contradicting upper semicontinuity of the function $f$ on the set
$\mathcal{A}_{k+1}$. $\square$ \medskip

An essential property of $\mu$\nobreakdash-\hspace{0pt}compact sets is presented in the
following proposition.\smallskip

\begin{property}\label{p-1} \emph{Let $\mathcal{A}$ be a convex $\mu$\nobreakdash-\hspace{0pt}compact
set and let $f$ be an upper semicontinuous upper bounded function on
a closed subset $\,\mathcal{B}\subset\mathcal{A}$. Then the function
\begin{equation}\label{f-mu}
\hat{f}^{\mu}_{\mathcal{B}}(x)=\sup_{\mu\in
M_{x}(\mathcal{B})}\int_{\mathcal{B}} f(y) \mu(dy)
\end{equation}
is upper semicontinuous and $\mu$\nobreakdash-\hspace{0pt}concave on the set
$\,\overline{\mathrm{co}}(\mathcal{B})$. For arbitrary
$x\in\overline{\mathrm{co}}(\mathcal{B})$ the supremum in the
definition of the value $\hat{f}^{\mu}_{\mathcal{B}}(x)$ is achieved
at a particular measure in $M_{x}(\mathcal{B})$.}
\end{property}
\medskip

This property provides generalization of several results well known
for compact convex sets to $\mu$\nobreakdash-\hspace{0pt}compact convex sets
\cite[Proposition 6, Corollary 2]{P&Sh}. It becomes no valid after
slight relaxing of the $\mu$\nobreakdash-\hspace{0pt}compactness assumption to pointwise
$\mu$\nobreakdash-\hspace{0pt}compactness \cite[Proposition 7]{P&Sh}.
The proof of Proposition \ref{p-1} is placed in the Appendix. \medskip

An another important technical tool is presented in the following
proposition. \smallskip

\begin{property}\label{p-2}
\emph{Let $\mathcal{A}$ be a convex $\mu$\nobreakdash-\hspace{0pt}compact \footnote{The
$\mu$\nobreakdash-\hspace{0pt}compactness assumption is used only to guarantee
$\mathbf{b}(M(\mathcal{B}))=\overline{\mathrm{co}}(\mathcal{B})$ by
means of Lemma \ref{bar-dec}.} set and let $f$ be a lower
semicontinuous lower bounded function on a closed subset
$\,\mathcal{B}\subseteq\mathcal{A}$.}
\begin{enumerate}[A)]
    \item \emph{If the map $M(\mathcal{B})\ni\mu\mapsto \mathbf{b}(\mu)\in
\overline{\mathrm{co}}(\mathcal{B})$ is open then the function $\hat{f}^{\mu}_{\mathcal{B}}$ defined by (\ref{f-mu})
is lower semicontinuous and $\mu$\nobreakdash-\hspace{0pt}concave on the set
$\;\overline{\mathrm{co}}(\mathcal{B})$.}
    \item \emph{If the map $M^a(\mathcal{B})\ni\mu\mapsto
\mathbf{b}(\mu)\in \sigma\textup{-}\mathrm{co}(\mathcal{B})$ is open
then the $\sigma$\nobreakdash-\hspace{0pt}concave  function
$$
\hat{f}^{\sigma}_{\mathcal{B}}(x)=\sup_{\mu\in
M^a_{x}(\mathcal{B})}\int_{\mathcal{B}} f(y)
\mu(dy)=\sup_{\{\pi_i,x_i\}\in M^a_{x}(\mathcal{B})}\sum_i \pi_i
f(x_i)
$$
is lower semicontinuous on the set
$\;\sigma\textup{-}\mathrm{co}(\mathcal{B})$. If, in addition,
$\;\sigma\textup{-}\mathrm{co}(\mathcal{B})=\overline{\mathrm{co}}(\mathcal{B})$
then the function $\hat{f}^{\sigma}_{\mathcal{B}}$ coincides with
the function $\hat{f}^{\mu}_{\mathcal{B}}$ defined by (\ref{f-mu}).}
\end{enumerate}
\end{property} \medskip

The proof of Proposition \ref{p-2} is placed in the Appendix. \medskip

\begin{remark}\label{p-1-2-r}
If $f$ is bounded and upper semicontinuous function on a closed subset $\mathcal{B}$ of a convex $\mu$\nobreakdash-\hspace{0pt}compact set $\mathcal{A}$ such that $\sigma\textup{-}\mathrm{co}(\mathcal{B})=\overline{\mathrm{co}}(\mathcal{B})$ then the above defined functions $\hat{f}^{\sigma}_{\mathcal{B}}$ and $\hat{f}^{\mu}_{\mathcal{B}}$ do not coincide in general (see the example in \cite[Remark 9]{Sh-9}). Thus the assertion of Propositions \ref{p-1} does not hold for the function $\hat{f}^{\sigma}_{\mathcal{B}}$ (since $\mu$\nobreakdash-\hspace{0pt}concavity of $\hat{f}^{\sigma}_{\mathcal{B}}$ implies $\hat{f}^{\sigma}_{\mathcal{B}}=\hat{f}^{\mu}_{\mathcal{B}}$).
\end{remark} \medskip

Propositions \ref{p-1} and \ref{p-2} have the obvious corollary.
\medskip

\begin{corollary}\label{p-1-2-c}
\emph{Let $\,\mathcal{B}$ be a closed subset of a convex $\mu$\nobreakdash-\hspace{0pt}compact
set $\mathcal{A}$.}
\begin{enumerate}[A)]
\item \emph{If $\mathcal{A}=\overline{\mathrm{co}}(\mathcal{B})$
and the map $M(\mathcal{B})\ni\mu\mapsto
\mathbf{b}(\mu)\in\mathcal{A}$ is open then
$\hat{f}^{\mu}_{\mathcal{B}}\in C(\mathcal{A})$ for any $f\in
C(\mathcal{B})$. }

\item \emph{If $\mathcal{A}=\sigma\textup{-}\mathrm{co}(\mathcal{B})$
and the map $M^a(\mathcal{B})\ni\mu\mapsto
\mathbf{b}(\mu)\in\mathcal{A}$ is open then
$\hat{f}^{\sigma}_{\mathcal{B}}=\hat{f}^{\mu}_{\mathcal{B}}\in
C(\mathcal{A})$ for any $f\in C(\mathcal{B})$.}
\end{enumerate}
\end{corollary}
\smallskip

If $\mathcal{A}$ is a stable convex $\mu$\nobreakdash-\hspace{0pt}compact set then the set
$\mathrm{extr}\mathcal{A}$ is closed and  the generalized
Vesterstrom-O'Brien theorem (\cite[Theorem 1]{P&Sh}) implies
openness of the surjective map
$M(\mathrm{extr}\mathcal{A})\ni\mu\mapsto\mathbf{b}(\mu)\in\mathcal{A}$,
hence Corollary \ref{p-1-2-c}A shows that an arbitrary function $f$
in $C(\mathrm{extr}\mathcal{A})$ has continuous bounded concave
extension $\hat{f}^{\mu}_{\mathrm{extr}\mathcal{A}}$ to the set
$\mathcal{A}$. This property does not hold in general for stable
convex sets, which are not $\mu$\nobreakdash-\hspace{0pt}compact (see Example 1 in
\cite{P&Sh}).\medskip

Corollary \ref{p-1-2-c}B plays an essential role in this paper due
to the following observation.\smallskip

\begin{property}\label{basic}
\textit{Let $\,\mathcal{A}_1$ be a closed subset of a stable convex
$\mu$\nobreakdash-\hspace{0pt}compact set $\,\mathcal{A}$ such that
$\,\mathcal{A}=\sigma\textup{-}\mathrm{co}(\mathcal{A}_1)$ and
$\,\{\mathcal{A}_{k}\}$ be the family of subsets defined by
(\ref{a-k}). If the map
\begin{equation}\label{ssp-1}
M^a(\mathcal{A}_{k})\ni\mu\mapsto\mathbf{b}(\mu)\in\mathcal{A}
\end{equation}
is open for $\,k=1$ then this map is open for all $\,k\in\mathbb{N}$.}
\end{property} \medskip

The proof of Proposition \ref{basic} is placed in the Appendix. \medskip

By the generalized Vesterstrom-O'Brien theorem stability of a convex $\mu$\nobreakdash-\hspace{0pt}compact
set $\mathcal{A}$ is equivalent to openness of the map
$M(\mathrm{extr}\mathcal{A})\ni\mu\mapsto\mathbf{b}(\mu)\in\mathcal{A}$.
By Lemma \ref{density} the last property implies openness of the map
$M^a(\mathrm{extr}\mathcal{A})\ni\mu\mapsto\mathbf{b}(\mu)\in\mathcal{A}$.
Hence we obtain from Proposition \ref{basic} the following
assertion. \medskip

\begin{corollary}\label{basic-c}
\textit{Let $\,\mathcal{A}$ be a stable convex
$\mu$\nobreakdash-\hspace{0pt}compact set such that
$\,\mathcal{A}=\sigma\textup{-}\mathrm{co}(\mathrm{extr}\mathcal{A})$
and $\{\mathcal{A}_{k}\}$ be the family of subsets defined by
(\ref{a-k}) with $\,\mathcal{A}_{1}=\mathrm{extr}\mathcal{A}$. Then
map (\ref{ssp-1}) is open for all $\,k\in\mathbb{N}$.}
\end{corollary}\medskip

\begin{remark}\label{basic-c-r}
If $\mathcal{A}$ is the stable convex $\mu$\nobreakdash-\hspace{0pt}compact set of quantum states (see Section 6) and
$\mathcal{A}_{1}=\mathrm{extr}\mathcal{A}$ then openness of the maps (\ref{ssp-1}) and
\begin{equation}\label{ssp-2}
M(\mathcal{A}_{k})\ni\mu\mapsto\mathbf{b}(\mu)\in\mathcal{A}
\end{equation}
are proved in \cite{Sh-11} by using special structure of this set and called strong stability property. By Corollary \ref{basic-c-r} to prove that
this strong stability property follows from stability it suffices to show that openness of map (\ref{ssp-2}) follows from openness of  map (\ref{ssp-1}). In \cite{Sh-11} this is made by proving density of the set $M^a_x(\mathcal{A}_{k})$ in $M_x(\mathcal{A}_{k})$ for all $x\in\mathcal{A}$.

\textbf{Question 1.}  Let $\mathcal{A}$ be a convex
$\mu$\nobreakdash-\hspace{0pt}compact set and $\{\mathcal{A}_{k}\}$ be
the family of subsets defined by (\ref{a-k}) with
$\mathcal{A}_{1}=\mathrm{extr}\mathcal{A}$. Does stability of the
set $\mathcal{A}$ imply openness of map (\ref{ssp-2}) for all $k$?
\end{remark}\smallskip

A positive answer on this question can be used to strengthen Theorem \ref{main+} in Section 5 (see Remark \ref{gcc-r} after this theorem).


\section{Special approximation of concave functions}

Throughout this section  we will assume that $f$ is a concave nonnegative function on a convex $\mu$\nobreakdash-\hspace{0pt}compact
set $\mathcal{A}$ having universally measurable restrictions
to subsets of the family $\{\mathcal{A}_{k}\}$ defined by (\ref{a-k}) with
$\mathcal{A}_{1}=\mathrm{cl}(\mathrm{extr}\mathcal{A})$. By Lemma
\ref{a-k-closed} this family consists of closed subsets. Possible generalizations are mentioned in Remark \ref{limit-r} at the end of this section.\smallskip

Since for arbitrary $x$ in $\mathcal{A}$ the set $M_{x}(\mathcal{A}_{1})$
is not empty by Proposition 5 in \cite{P&Sh}, for given natural $k$
we can consider the concave nonnegative function
\begin{equation}\label{f-mu-k}
\mathcal{A}\ni x\mapsto\hat{f}_{k}^{\mu}(x)=\sup_{\mu\in
 M_{x}(\mathcal{A}_{k})} \int_{\mathcal{A}_{k}}f(y)\mu(dy).
\end{equation}
If the function $f$ is $\mu$\nobreakdash-\hspace{0pt}concave on the
set $\mathcal{A}$ then
\begin{equation}\label{one}
\hat{f}_{k}^{\mu}\leq f\quad
\textup{and}\quad\hat{f}_{k}^{\mu}|_{\mathcal{A}_{k}}=f|_{\mathcal{A}_{k}},
\end{equation}
hence the function $\hat{f}_{k}^{\mu}$ can be considered as a
concave extension of the function $f|_{\mathcal{A}_{k}}$ to the set
$\mathcal{A}$. If the function $f$ has upper semicontinuous bounded
restriction to the set $\mathcal{A}_{k}$ then the bounded function
$\hat{f}_{k}^{\mu}$ is upper semicontinuous and
$\mu$\nobreakdash-\hspace{0pt}concave on the set $\mathcal{A}$ by
Proposition \ref{p-1}. Hence in this case the function
$\hat{f}_{k}^{\mu}$ is the minimal
$\mu$\nobreakdash-\hspace{0pt}concave extension of the function
$f|_{\mathcal{A}_{k}}$ to the set $\mathcal{A}$.

The sequence  $\,\{\hat{f}_{k}^{\mu}\}$ is nondecreasing and its pointwise limit $\,\hat{f}_{\ast}^{\mu}\doteq\sup_k\hat{f}_{k}^{\mu}$ is a concave function on $\mathcal{A}$. If the function $f$ is $\mu$\nobreakdash-\hspace{0pt}concave then (\ref{one}) implies
\begin{equation}\label{one+}
\hat{f}_{\ast}^{\mu}\leq f\quad
\textup{and}\quad f_{\ast}^{\mu}|_{\mathcal{A}_{\ast}}=f|_{\mathcal{A}_{\ast}},\quad
\textup{where}\quad
\mathcal{A}_{\ast}=\bigcup_{k=1}^{+\infty}\mathcal{A}_{k}.
\end{equation}

\textbf{Question 2.} Under what conditions do the functions  $\,\hat{f}_{\ast}^{\mu}$ and  $\,f$ coincide? \medskip

A partial answer on this question can be obtained in the case
$\mathcal{A}=\sigma\textup{-}\mathrm{co}(\mathcal{A}_1)$.\smallskip

In this case for given
natural $k$  one can consider the $\sigma$\nobreakdash-\hspace{0pt}concave
nonnegative function
\begin{equation}\label{f-sigma-k}
\mathcal{A}\ni
x\mapsto\hat{f}_{k}^{\sigma}(x)=\sup_{\{\pi_{i},x_{i}\}\in
M^{\mathrm{a}}_{x}(\mathcal{A}_{k})} \sum_{i}\pi_{i}f(x_{i}).
\end{equation}
By the construction $\,\hat{f}_{k}^{\sigma}\leq\hat{f}_{k}^{\mu}\,$. Since the function $f$ is $\sigma$\nobreakdash-\hspace{0pt}concave on
the set $\mathcal{A}$ by Lemma \ref{yensen}, we have
\begin{equation}\label{two}
\hat{f}_{k}^{\sigma}\leq f\quad
\textup{and}\quad\hat{f}_{k}^{\sigma}|_{\mathcal{A}_{k}}=f|_{\mathcal{A}_{k}}.
\end{equation}
Hence the function $\hat{f}_{k}^{\sigma}$ is the
minimal $\sigma$\nobreakdash-\hspace{0pt}concave extension of the
function $f|_{\mathcal{A}_{k}}$ to the set $\mathcal{A}$.\medskip

The sequence
$\,\{\hat{f}_{k}^{\sigma}\}$ is nondecreasing and its pointwise limit $\,\hat{f}_{\ast}^{\sigma}\doteq\sup_k\hat{f}_{k}^{\sigma}$ is a concave function on $\mathcal{A}$ such that $\hat{f}_{\ast}^{\sigma}\leq\hat{f}_{\ast}^{\mu}$. It follows from (\ref{two}) that relations (\ref{one+}) hold with $\hat{f}_{\ast}^{\sigma}$ instead of $\hat{f}_{\ast}^{\mu}$. \medskip

\begin{property}\label{limit} \emph{If $\,\mathcal{A}=\sigma\textup{-}\mathrm{co}(\mathcal{A}_1)$ and the function $\,f$ is lower
semicontinuous then
$$\,\hat{f}_{\ast}^{\mu}=\hat{f}_{\ast}^{\sigma}=f.$$}
\end{property}
\smallskip

\textbf{Proof.} By Lemma \ref{yensen} the function $f$ is
$\mu$\nobreakdash-\hspace{0pt}concave. Hence (\ref{one})  holds for
all $k$ and to prove
$\hat{f}_{\ast}^{\sigma}=\hat{f}_{\ast}^{\mu}=f$ it is sufficient to
show that $\hat{f}_{\ast}^{\sigma}=f$.

Let $x_{0}$ be an arbitrary point in $\mathcal{A}$. Then
$x_{0}=\sum_{i=1}^{+\infty}\pi_i y_i$, where
$\{\pi_i\}\in\mathfrak{P}_{+\infty}$ and $\{y_i\}\in\mathcal{A}_1$.
Let $x_{n}=(\lambda_n)^{-1}\sum_{i=1}^{n}\pi_i y_i$ and
$y_{n}=(1-\lambda_n)^{-1}\sum_{i>n}\pi_i y_i$, where
$\lambda_n=\sum_{i=1}^{n}\pi_i$. The sequence $\{x_{n}\}$ belongs to
the set $\mathcal{A}_{\ast}$ and converges to the point $x_{0}$.

For each $n$ we have $x_0=\lambda_n x_n + (1-\lambda_n)y_n$ and
hence
$\hat{f}_{\ast}^{\sigma}(x_{0})\geq\lambda_{n}\hat{f}_{\ast}^{\sigma}(x_{n})=\lambda_{n}f(x_{n})$
by concavity and nonnegativity of the function
$\hat{f}_{\ast}^{\sigma}$. This implies
$\limsup_{n\rightarrow+\infty}f(x_{n})\leq
\hat{f}_{\ast}^{\sigma}(x_{0})$. By lower semicontinuity of the
function $f$ we have $f(x_{0})\leq\hat{f}_{\ast}^{\sigma}(x_{0})$
and hence $f(x_{0})=\hat{f}_{\ast}^{\sigma}(x_{0})$. $\square$

Lemma \ref{a-k-closed}B, Proposition \ref{p-1}, Corollary
\ref{p-1-2-c}B and Corollary \ref{basic-c} imply the following
observation, providing usefulness of the approximating sequences $\{\hat{f}_{k}^{\mu}\}$ and $\{\hat{f}_{k}^{\sigma}\}$ for our purposes.\smallskip

\begin{property}\label{limit+} \emph{If the function $f$ has continuous
restriction to the set $\,\mathcal{A}_{k}$ for each $k$ then the
function $\hat{f}_{k}^{\mu}$ is bounded and upper semicontinuous for
each $\,k$. If, in addition, the set $\,\mathcal{A}$ is stable and $\,\mathcal{A}=\sigma\textup{-}\mathrm{co}(\mathcal{A}_1)$ then
$\,\hat{f}_{k}^{\sigma}=\hat{f}_{k}^{\mu}\in C(\mathcal{A})$ for each
$\,k$}.
\end{property}
\smallskip

\begin{remark}\label{limit-r}
The above constructions can be generalized by considering the family
$\{\mathcal{A}_{k}\}$ produced by an arbitrary closed subset
$\mathcal{A}_1$ of $\mathcal{A}$ such that
$\mathcal{A}=\overline{\mathrm{co}}(\mathcal{A}_1)$. The all results
remain valid in this case excepting the second assertion of
Proposition \ref{limit+}, in which the requirement of openness of
the map
$M^a(\mathcal{A}_{1})\ni\mu\mapsto\mathbf{b}(\mu)\in\mathcal{A}$
must be added. This can be shown by applying Proposition \ref{basic}
instead of Corollary \ref{basic-c}.
\end{remark}

\section{Continuity conditions}

Let $f$ be a concave nonnegative function on a convex $\mu$\nobreakdash-\hspace{0pt}compact
set $\mathcal{A}$. In this section we consider conditions of
continuity of this function on subsets of $\mathcal{A}$ assuming
that there exists a closed subset $\mathcal{A}_{1}\subset\mathcal{A}$
such that $\mathcal{A}=\overline{\mathrm{co}}(\mathcal{A}_1)$ and
\begin{equation}\label{b-a}
f|_{\mathcal{A}_{k}}\;\,\textrm{is}\;\, \textrm{continuous}\;\,
\textrm{for}\;\, \textrm{each}\;\, \textrm{natural}\;\,  k,
\end{equation}
where $\mathcal{A}_{k}$ is the subset of $\mathcal{A}$ defined by
(\ref{a-k}). This assumption with
$\mathcal{A}_{1}=\mathrm{cl}(\mathrm{extr}\mathcal{A})$ has a
physical motivation (see Section 6). Sometimes it can be
reduced to continuity and boundedness of $f|_{\mathcal{A}_{1}}$ (see the proof of
Lemma \ref{l-lemma} below).


\subsection{The case $\mathcal{A}=\sigma\textup{-}\mathrm{co}(\mathrm{cl}(\mathrm{extr}\mathcal{A}))$}

The results of the previous sections imply the following continuity
condition.
\medskip

\begin{theorem}\label{main}
\emph{Let $\mathcal{A}$ be a convex $\mu$\nobreakdash-\hspace{0pt}compact set such that
$\mathcal{A}=\sigma\textup{-}\mathrm{co}(\mathrm{cl}(\mathrm{extr}\mathcal{A}))$.
Let $f$ be a concave nonnegative function on the set $\mathcal{A}$ such that
assumption (\ref{b-a}) holds with
$\mathcal{A}_{1}=\mathrm{cl}(\mathrm{extr}\mathcal{A})$.  Assume that one of
the following conditions is valid:}

\begin{enumerate}[a)]
    \item \emph{the set $\mathcal{A}$ is stable,}
    \item \emph{the function $f$ is lower semicontinuous.}
\end{enumerate}
\emph{Then the function $f$ is continuous
on a subset $\,\mathcal{B}\subseteq\mathcal{A}$ if
\begin{equation}\label{scc}
\lim_{k\rightarrow+\infty}\sup_{x\in\mathcal{B}}\Delta^\sigma_k(x|f)=0,
\quad where \quad \Delta^\sigma_k(x|f)=\inf_{\{\pi_i, x_i\}\in
M^a_x(\mathcal{A}_k)}\left[f(x)-\sum_i \pi_i f(x_i)\right].
\end{equation}}
\emph{If the both above conditions $\,\mathrm{a)}$ and
$\,\mathrm{b)}$ are valid then (\ref{scc}) is a necessary and sufficient condition
of continuity of the function $f$ on a compact subset
$\,\mathcal{B}\subset\mathcal{A}$.}
\end{theorem}
\medskip

\begin{remark}\label{scc-r}
Since $\Delta^\sigma_k(x|f)=f-\hat{f}^\sigma_k$ and $\hat{f}^\sigma_k\leq\hat{f}^\mu_k$, where $\hat{f}^\mu_k$ and $\hat{f}^\sigma_k$ are functions defined by (\ref{f-mu-k}) and (\ref{f-sigma-k}), condition
(\ref{scc}) means uniform convergence of the sequences
$\{\hat{f}^\sigma_k\}$ and $\{\hat{f}^\mu_k\}$ to the function $f$ on the subset $\mathcal{B}$.
\end{remark}
\smallskip

\begin{remark}\label{scc-r+}
Applications of the above continuity condition are based on
possibility to find for a given concave function $f$ a suitable
upper bound for the value in the square brackets in (\ref{scc}) (see Example \ref{main-e}  below and
Section 6).
\end{remark}
\smallskip

\textbf{Proof.} If the set $\mathcal{A}$ is stable then
$\hat{f}^\mu_k=\hat{f}^\sigma_k\in C(\mathcal{A})$ for all $k$ by
Proposition \ref{limit+}. By Remark \ref{scc-r} condition
(\ref{scc}) implies continuity of the function $f$ on the subset
$\mathcal{B}$.

If the function $f$ is lower semicontinuous then continuity of the
function $f$ on the subset $\mathcal{B}$ can be verified by showing
its upper semicontinuity and boundedness on this set. By Remark
\ref{scc-r} the last property follows from condition (\ref{scc})
since by Proposition \ref{limit+} the sequence $\{\hat{f}^\mu_k\}$ consists of upper
semicontinuous bounded functions.

By Propositions \ref{limit} and \ref{limit+} the last
assertion of the theorem follows from Dini's lemma and Remark \ref{scc-r}. $\square$\smallskip

\begin{example}\label{main-e}
The Shannon entropy is a concave lower semicontinuous function on
the set
$\mathfrak{P}_{+\infty}=\left\{\,x=\{x^j\}_{j=1}^{+\infty}\in\ell_1\,|\,x^j\geq\;
0,\; \forall j,\;\, \sum_{j=1}^{+\infty}x^j=1\,\right\}$ of all countable
probability distributions defined as follows
$$
S\left(\{x^j\}_{j=1}^{+\infty}\right)=-\sum_{j=1}^{+\infty}x^j\ln
x^j.
$$
This function is nonnegative and takes the value  $+\infty$ on a
dense subset of $\mathfrak{P}_{+\infty}$.

As mentioned in Section 2 the convex set
$\mathfrak{P}_{+\infty}$ is stable and $\mu$-compact. The set
$\mathrm{extr}\,\mathfrak{P}_{+\infty}$ consists of "degenerate"
distributions having $"1"$ at some position and $"0"$ on other
places. It is clear that
$\mathfrak{P}_{+\infty}=\sigma\textup{-}\mathrm{co}\left(\mathrm{extr}\,\mathfrak{P}_{+\infty}\right)$
and that the function $x\mapsto S(x)$ has continuous restriction to the set
$$
(\mathfrak{P}_{+\infty})_k=\left\{\left.\sum_{i=1}^{k}\pi_i x_i\,
\right| \{\pi_i\}\in\mathfrak{P}_k,\,
\{x_i\}\subset\mathrm{extr}\,\mathfrak{P}_{+\infty}\right\}
$$
for each $\,k\in\mathbb{N}$. If $f=S$ then the value in the squire
brackets in (\ref{scc}) can be expressed as follows
$$
S(x)-\sum_i \pi_i S(x_i)=\sum_{i}\pi_{i}S(x_{i}\|\,x),
$$
where $S(\cdot\|\cdot)$ is the relative entropy (Kullback-Leibler
distance \cite{Kullback}) defined for arbitrary distributions
$x=\{x^j\}_{j=1}^{+\infty}$ and $y=\{y^j\}_{j=1}^{+\infty}$ in
$\mathfrak{P}_{+\infty}$ by the formula
$$
S(x\|\,y)=\left\{\begin{array}{cc}
  \sum_{i=1}^{+\infty}x^j \ln
(x^j/y^j), & \{y^j=0\}\Rightarrow\{x^j=0\} \\
  +\infty, & \textrm{otherwise}
\end{array}\right..
$$

Thus Theorem \ref{main} implies the following continuity condition
for the Shannon entropy.\smallskip

\emph{The function $\,x\mapsto S(x)$ is continuous on a compact
subset $\,\mathfrak{P}\subseteq\mathfrak{P}_{+\infty}$ if and only
if
\begin{equation}\label{S-cont-cond}
\lim_{k\rightarrow+\infty}\sup_{x\in\mathfrak{P}}\Delta^{\sigma}_k(x|S)=0,
\quad where \quad \Delta^{\sigma}_k(x|S)=\inf_{\{\pi_i, x_i\}\in
M^a_x((\mathfrak{P}_{+\infty})_k)}\sum_{i}\pi_{i}S(x_{i}\|\,x).
\end{equation}}

This condition can be applied directly by using well
studied properties of the relative entropy. For example, by joint
convexity and lower semicontinuity of the relative entropy validity of (\ref{S-cont-cond})
for convex subsets $\mathfrak{P}'$ and $\mathfrak{P}''$ of
$\mathfrak{P}_{+\infty}$ implies validity of (\ref{S-cont-cond}) for
their convex closure
$\overline{\mathrm{co}}(\mathfrak{P}'\cup\mathfrak{P}'')$. Hence, we
can conclude that \emph{continuity of the Shannon entropy on convex
closed subsets $\,\mathfrak{P}'$ and $\,\mathfrak{P}''$ implies its
continuity on their convex closure
$\,\overline{\mathrm{co}}(\mathfrak{P}'\cup\mathfrak{P}'')$}.\footnote{It
is possible to show that continuity of the Shannon entropy on a
convex subset of $\mathfrak{P}_{+\infty}$ implies relative compactness
of this subset.}\smallskip

The above continuity condition can be also applied by using the
estimation
\begin{equation}\label{d-est}
\Delta^{\sigma}_k(x|S)\leq S(k(x)),\quad k\in\mathbb{N},
\end{equation}
where $\,k(x)$ is a distribution obtained by $k$-order
coarse-graining from the distribution $x$, that is
$\,(k(x))^j=x^{(j-1)k+1}+...+x^{jk}\,$ for all $\,j=1,2,...$. This
estimation is proved by using the decomposition
$x=\sum_{i=1}^{+\infty}\lambda^k_{i}p^k_{i}(x)$, where
$\lambda^k_{i}=(k(x))^i$ and $p^k_{i}(x)$ is a distribution such
that $(p^k_{i}(x))^{j}=(\lambda^k_{i})^{-1}x^{j}$ for
$j=\overline{(i-1)k+1, ik}\,$ and $(p^k_{i}(x))^{j}=0$ for others $j$,
since it is easy to verify that
$\sum_{i=1}^{+\infty}\lambda^k_{i}S(p^k_{i}(x)\|\,x)=\sum_{i=1}^{+\infty}\lambda^k_{i}(-\ln\lambda^k_{i})=S(k(x))$.\medskip

The above continuity condition  and estimation (\ref{d-est}) imply
the following assertion.\smallskip

\emph{Let $\,x_0$ be a distribution in $\,\mathfrak{P}_{+\infty}$
with finite Shannon entropy, then the Shannon entropy is
continuous on the set
\begin{equation}\label{m-set}
\left\{\,x\in\mathfrak{P}_{+\infty}\,|\;x\prec x_0\,\right\},
\end{equation}
where $x\prec y$ means that the distribution
$\,y=\{y^j\}_{j=1}^{+\infty}$ is more chaotic than the distribution
$\,x=\{x^j\}_{j=1}^{+\infty}$ in the Uhlmann sense \cite{U+,W+}, that
is $\,\sum_{j=1}^{n}x^{j}\geq\sum_{j=1}^{n}y^{j}$ for each
natural $n$ provided the sequences $\,\{x^j\}_{j=1}^{+\infty}$ and
$\,\{y^j\}_{j=1}^{+\infty}$ are arranged in nonincreasing
order.}\footnote{The order $"\prec"$ is converse to the majorization
order used in linear algebra \cite{Bhatia}.}\smallskip

Indeed, assuming that the elements of $x$ and $x_0$ are arranged in nonincreasing
order we have $x\prec x_0\Rightarrow k(x)\prec k(x_0)\Rightarrow S(k(x))\leq S(k(x_0))$ by Shur concavity of the Shannon entropy
\cite{W+}. Hence validity of (\ref{S-cont-cond}) for set
(\ref{m-set}) follows from
(\ref{d-est}) and the easily verified implication
$\,S(x_0)<+\infty\,\Rightarrow\, \lim_{k\rightarrow+\infty} S(k(x_0))=0$.
\end{example}

\subsection{Possible generalizations}

Note first that Theorem \ref{main} can be generalized by replacing
the family $\{\mathcal{A}_k\}$ produced by the set
$\mathcal{A}_1=\mathrm{cl}(\mathrm{extr}(\mathcal{A}))$ by a family
$\{\mathcal{A}_k\}$ produced by an arbitrary closed subset
$\mathcal{A}_1$ of $\mathcal{A}$ such that
$\mathcal{A}=\sigma\textup{-}\mathrm{co}(\mathcal{A}_1)$. By Remark
\ref{limit-r} the only necessary modification of Theorem \ref{main}
under this replacement consists in the additional requirement of
openness of the map
$M^a(\mathcal{A}_{1})\ni\mu\mapsto\mathbf{b}(\mu)\in\mathcal{A}\,$
in condition $\mathrm{a})$.
\medskip

Without the assumption
$\mathcal{A}=\sigma\textup{-}\mathrm{co}(\mathcal{A}_1)$ the
following continuity condition can be proved.\medskip

\begin{theorem}\label{main+}
\emph{Let $\mathcal{A}_{1}$ be a closed subset of a convex
$\mu$\nobreakdash-\hspace{0pt}compact set $\mathcal{A}$ such that
$\mathcal{A}=\overline{\mathrm{co}}(\mathcal{A}_1)$ (in particular,
$\mathcal{A}_1=\mathrm{cl}(\mathrm{extr}(\mathcal{A}))$). Let $f$ be
a concave lower semicontinuous nonnegative function on the set $\mathcal{A}$
satisfying assumption (\ref{b-a}). Then the function $f$ is
continuous on a subset $\,\mathcal{B}\subseteq\mathcal{A}$ if
\begin{equation}\label{gcc}
\lim_{k\rightarrow+\infty}\sup_{x\in\mathcal{B}}\Delta^\mu_k(x|f)=0,
\quad where \quad \Delta^\mu_k(x|f)=\inf_{\mu\in
M_x(\mathcal{A}_k)}\left[f(x)-\int_{\mathcal{A}_k}f(y)\mu(dy)\right].
\end{equation}
Condition (\ref{gcc}) can be replaced by the following one
\begin{equation}\label{gcc+}
\lim_{k\rightarrow+\infty}\sup_{x\in\mathcal{B}_0}\Delta^\sigma_k(x|f)=0,
\quad where \quad \Delta^\sigma_k(x|f)=\inf_{\{\pi_i, x_i\}\in
M^a_x(\mathcal{A}_k)}\left[f(x)-\sum_i \pi_i f(x_i)\right]
\end{equation}
and $\,\mathcal{B}_0$ is an arbitrary subset of
$\,\sigma\textup{-}\mathrm{co}(\mathcal{A}_1)$ such that
$\,\mathcal{B}\subseteq\mathrm{cl}(\mathcal{B}_0)$.}\smallskip
\end{theorem}\medskip

\textbf{Proof.} By Proposition \ref{limit+} the function
$\hat{f}^\mu_k$ defined by formula (\ref{f-mu-k}) is upper semicontinuous and bounded
for each $k$. Since $\Delta^\mu_k(x|f)=f-\hat{f}^\mu_k$,  condition
(\ref{gcc}) means uniform convergence of the sequence
$\{\hat{f}^\mu_k\}$ to the function $f$ on the subset $\mathcal{B}$,
which implies upper semicontinuity and boundedness of the lower semicontinuous function
$f$ on the subset $\mathcal{B}$.

Let $\mathcal{B}_0$ be a subset of
$\sigma\textup{-}\mathrm{co}(\mathcal{A}_1)$ such that
$\mathcal{B}\subseteq\mathrm{cl}(\mathcal{B}_0)$. On this subset the
function $\hat{f}^\sigma_k$ is well defined by formula (\ref{f-sigma-k}) for each $k$. Since $\Delta^\sigma_k(x|f)=f-\hat{f}^\sigma_k$ and $\hat{f}^\sigma_k\leq\hat{f}^\mu_k$, condition (\ref{gcc+})  guarantees uniform convergence
of the sequence $\{\hat{f}^\mu_k\}$ to the function $f$ on the
subset $\mathcal{B}_0$, which implies uniform
convergence of the sequence $\{\hat{f}^\mu_k\}$ to the function $f$
on the subset $\mathrm{cl}(\mathcal{B}_0)$, since the function
$\Delta^\mu_k(x|f)=f-\hat{f}^\mu_k$ is lower semicontinuous (as a
difference between lower semicontinuous and bounded upper
semicontinuous functions). $\square$ \smallskip

\begin{remark}\label{gcc-r}
If the set $\mathcal{A}$ is stable, $\mathcal{A}_{1}=\mathrm{extr}(\mathcal{A})$ and positive answers on the above Questions 1 and 2 (stated respectively in Sections 3 and 4) hold then (\ref{gcc}) is a necessary and sufficient condition
of continuity of the function $f$ on a compact subset
$\,\mathcal{B}\subset\mathcal{A}$. Necessity of condition (\ref{gcc}) in this case can shown by using Corollary
\ref{p-1-2-c}A and Dini's lemma.
\end{remark}\medskip

Theorem \ref{main+} can be applied to analysis of concave functions on the stable convex $\mu$\nobreakdash-\hspace{0pt}compact set of probability measures on a complete separable metric space having continuous restrictions to the subset of measures supported by $\leq k$ atoms for all $k$.

\section{Applications in quantum physics}

The notion of a quantum state plays a central role in the
statistical structure of quantum theory \cite{H-SSQT}. In this
section we consider applications of the continuity conditions
obtained in the previous section to analysis of local continuity of
several entropic characteristics -- the particular concave functions
on the convex set of all quantum states.\smallskip

Let $\mathcal{H}$ be a separable Hilbert space,
$\mathfrak{T}( \mathcal{H})$ -- the Banach space of all
trace\nobreakdash-\hspace{0pt}class operators in $\mathcal{H}$ with
the trace norm, containing the cone
$\mathfrak{T}_{+}(\mathcal{H})$ of all positive
trace\nobreakdash-\hspace{0pt}class operators.

The closed convex set
$$
\mathfrak{S}(\mathcal{H})=\{A\in\mathfrak{T}_{+}(\mathcal{H})\,|\,\mathrm{Tr}A=1\}
$$
is a complete separable metric space with the metric defined by the
trace norm. Operators in $\mathfrak{S}(\mathcal{H})$ are denoted
$\rho,\sigma,\omega,...$ and called density operators or \emph{quantum states}
since each density operator corresponds to a normal state on the algebra  of all bounded
operators \cite{B&R}.

It is essential that the convex set $\mathfrak{S}(\mathcal{H})$ is
stable and $\mu$\nobreakdash-\hspace{0pt}compact \cite{Sh-9} (the set
$\mathfrak{S}(\mathcal{H})$ is compact if and only if
$\dim\mathcal{H}<+\infty$). The set
$\mathrm{extr}\mathfrak{S}(\mathcal{H})$ of its extreme points consists
of one dimensional projectors -- \emph{pure} states. A pure state corresponding to a unit vector $|\varphi\rangle\in\mathcal{H}$ will be denoted $|\varphi\rangle\langle\varphi|$. By the spectral theorem an arbitrary  state $\rho$ can be represented
as follows $\rho=\sum_i \lambda_i |\varphi_i\rangle\langle \varphi_i|$, where
$\{|\varphi_i\rangle\}$ is the orthonormal basis of
eigenvectors  of the operator $\rho$ and
$\{\lambda_i\}$ is the corresponding sequence of
eigenvalues. Hence $\mathfrak{S}(\mathcal{H})=\sigma\textup{-}\mathrm{co}(\mathrm{extr}\mathfrak{S}(\mathcal{H}))$.

Rapid development of quantum information theory leads to discovery
of a whole number of important entropic and informational
characteristics of quantum systems, see e.g. \cite{H-SSQT,N&Ch}. Many of
them can be considered as functions on the set of quantum states. In
the finite dimensional case ($\dim\mathcal{H}<+\infty$) these
functions are generally bounded and continuous on the whole set of quantum
states, but in infinite dimensions their analytical properties are not so good.
For example, the von Neumann entropy is
a continuous bounded function on the set of quantum states of finite
dimensional quantum system, but it is discontinuous and takes the value $+\infty$ on a
dense subset of the set of all infinite dimensional quantum
states.\footnote{Moreover, the set of states with finite von Neumann entropy is
a first category subset of the set of all quantum states \cite{W}.}

Discontinuity and unboundedness of entropic characteristics lead to
technical problems in analysis of infinite dimensional quantum
systems. Moreover, they produce a number of "nonphysical" effects
such as infinite values of different capacities of a quantum channel and their
discontinuity as functions of a channel \cite{L&S,Sh-3}. But these difficulties can be partially overcome by using
local continuity conditions for entropic characteristics \cite{L,Sh-3,W}. For example, continuity of the von Neumann entropy
on the set of states of the system of quantum
oscillators with bounded mean energy provides many results concerning different characteristics of this system (see \cite{E&W} and references therein). Thus, study of local continuity properties of entropic
characteristics of quantum states is important for
rigorous analysis of infinite dimensional quantum systems.

Since $\mathfrak{S}(\mathcal{H})=\sigma\textup{-}\mathrm{co}(\mathrm{extr}\mathfrak{S}(\mathcal{H}))$ is a convex stable $\mu$\nobreakdash-\hspace{0pt}compact set,
we can apply the results of the previous sections to
study concave nonnegative functions on the set
$\mathfrak{S}(\mathcal{H})$ having  restrictions to the set
\begin{equation}\label{s-k}
\mathfrak{S}_k(\mathcal{H})=\left\{\left.\sum_{i=1}^{k}\pi_i
\rho_i\, \right| \{\pi_i\}\in\mathfrak{P}_k,\,
\{\rho_i\}\subset\mathrm{extr}\mathfrak{S}(\mathcal{H})\right\}
\end{equation}
with appropriate analytical properties for all $k$. Note that
$\mathfrak{S}_k(\mathcal{H})$ is the set of all quantum states having
rank $\leq k$ (as operators in $\mathcal{H}$), it can be considered as an union of all unitary
translations of the set $\mathfrak{S}(\mathcal{H}_k)$, where
$\mathcal{H}_k$ is a particular $k$-dimensional subspace of
$\mathcal{H}$.\medskip

Let $f$ be a concave nonnegative function on the set
$\mathfrak{S}(\mathcal{H})$. For given natural $k$ consider the
concave functions
$$
\hat{f}_{k}^{\mu}(\rho)=\sup_{\mu\in
M_{\rho}(\mathfrak{S}_k(\mathcal{H}))}
\int_{\mathfrak{S}_k(\mathcal{H})}f(\sigma)\mu(d\sigma)\quad\textrm{and}\quad\hat{f}_{k}^{\sigma}(\rho)=\sup_{\{\pi_{i},\rho_{i}\}\in
M^{a}_{\rho}(\mathfrak{S}_k(\mathcal{H}))} \sum_{i}\pi_{i}f(\rho_{i})
$$
on the set $\mathfrak{S}(\mathcal{H})$ (assuming that $f$ has
universally measurable restriction to the set
$\mathfrak{S}_k(\mathcal{H})$).\medskip

It is clear that $\,\hat{f}_{k}^{\sigma}\leq \hat{f}_{k}^{\mu}\,$.
Since the function $f$ is $\sigma$\nobreakdash-\hspace{0pt}concave by Lemma \ref{yensen}A we
have $\,\hat{f}_{k}^{\sigma}\leq f\,$ and
$\;\hat{f}_{k}^{\sigma}|_{\mathfrak{S}_k(\mathcal{H})}=f|_{\mathfrak{S}_k(\mathcal{H})}\,$.
If the function $f$ is $\mu$\nobreakdash-\hspace{0pt}concave (see conditions in Lemma
\ref{yensen}B) then $\,\hat{f}_{k}^{\mu}\leq f\,$ and
$\;\hat{f}_{k}^{\mu}|_{\mathfrak{S}_k(\mathcal{H})}=f|_{\mathfrak{S}_k(\mathcal{H})}\,$.
\medskip

The results of Sections 3 and 4 imply the following
observations.\medskip

\begin{property}\label{approx} \emph{Let $f$ be a concave  nonnegative function
on the set $\,\mathfrak{S}(\mathcal{H})$, taking finite value at least at one state.}
\begin{enumerate}[A)]
    \item \emph{If $f|_{\mathfrak{S}_k(\mathcal{H})}$ is upper
    semicontinuous for each $k$ then the function $\hat{f}_{k}^{\mu}$ is upper
    semicontinuous and bounded on the set $\,\mathfrak{S}(\mathcal{H})$ for each $k$.}
    \item \emph{If $f|_{\mathfrak{S}_k(\mathcal{H})}$ is lower
    semicontinuous for each $k$ then $\hat{f}_{k}^{\mu}=\hat{f}_{k}^{\sigma}$ and this function is lower
    semicontinuous on the set $\,\mathfrak{S}(\mathcal{H})$ for each $k$.}
    \item \emph{If $f|_{\mathfrak{S}_k(\mathcal{H})}$ is continuous for each $k$ then
    $\hat{f}_{k}^{\mu}=\hat{f}_{k}^{\sigma}\in C(\mathfrak{S}(\mathcal{H}))$ for each
    $k$.}
\end{enumerate}

\emph{If the function $f$ is lower semicontinuous on the set $\,\mathfrak{S}(\mathcal{H})$ then the nondecreasing sequence
$\{\hat{f}_{k}^{\mu}=\hat{f}_{k}^{\sigma}\}$ pointwise converges to
the function $f$.}\medskip
\end{property}

By Proposition \ref{approx} an arbitrary concave lower
semicontinuous nonnegative function $f$ on the set
$\mathfrak{S}(\mathcal{H})$ having continuous restriction to the
set $\mathfrak{S}_{k}(\mathcal{H})$ for each $k$ can be approximated
by the increasing sequence
of concave
continuous nonnegative bounded functions\break $f_k\doteq\hat{f}_{k}^{\mu}=\hat{f}_{k}^{\sigma}\,$ such that
$f_{k}|_{\mathfrak{S}_k(\mathcal{H})}=f|_{\mathfrak{S}_k(\mathcal{H})}$
for each $k$. Advantages of this approximation and its possible
applications are considered in \cite[Section 4]{Sh-11} and in
\cite[Section 6.2]{Sh-9}.
\medskip

Theorem \ref{main} implies the following continuity condition (extending the results of \cite{Sh-11}).

\medskip

\begin{property}\label{main-p}
\emph{Let $f$ be a concave nonnegative function on the set
$\,\mathfrak{S}(\mathcal{H})$ having continuous restriction to the
set $\,\mathfrak{S}_{k}(\mathcal{H})$ defined by (\ref{s-k}) for
each $k$. Then the function $f$ is continuous on a subset
$\,\mathfrak{S}\subseteq\mathfrak{S}(\mathcal{H})$ if
\begin{equation}\label{scc+}
\lim_{k\rightarrow+\infty}\sup_{\rho\in\mathfrak{S}}\Delta^\sigma_k(\rho|f)=0,
\quad where \quad \Delta^\sigma_k(\rho|f)=\inf_{\{\pi_i, \rho_i\}\in
M^a_\rho(\mathfrak{S}_{k}(\mathcal{H}))}\left[f(\rho)-\sum_i \pi_i
f(\rho_i)\right].
\end{equation}}
\emph{If the function $f$ is lower semicontinuous then (\ref{scc+})
is a necessary and sufficient condition of continuity of the
function $f$ on a compact subset
$\,\mathfrak{S}\subset\mathfrak{S}(\mathcal{H})$.}
\end{property}\medskip

The conditions of Proposition \ref{main-p} are valid for the
following well known characteristics of quantum states -- concave
lower semicontinuous nonnegative functions on the set
$\mathfrak{S}(\mathcal{H})$:
\begin{itemize}
    \item the quantum Renyi entropy
    $R_p(\rho)=\ln\mathrm{Tr}\rho^p/(1-p)$ of order
    $p\in(0,1]$ (the case $p=1$ corresponds to the von Neumann entropy $H(\rho)=-\mathrm{Tr}\rho\ln\rho$);
    \item the quantum mutual information $I(\rho, \Phi)$ of a quantum
    channel $\Phi$ (defined in Section 6.2);
    \item the output quantum Renyi entropy  $R_p(\Phi(\rho))$ of order
    $p\in(0,1]$ (in particular, the output von Neumann entropy $H(\Phi(\rho)))$
    of a quantum channel $\Phi$ satisfying the particular condition (see Section 6.3).
    \item the $\chi$-function (the constrained Holevo capacity) $\chi_{\Phi}(\rho)$ a quantum channel $\Phi$ satisfying the particular condition (see Section 6.3).
\end{itemize}

Below we consider applications of Proposition \ref{main-p} to the
above functions, reducing attention to the von Neumann entropy --
the most important version of the quantum Renyi entropy.

\subsection{The von
Neumann entropy}

Continuity conditions for the von Neumann entropy on subsets of $\mathfrak{T}_{+}(\mathcal{H})$ based on the above approximation technic are presented in \cite{Sh-11}. Here we consider the case of the von Neumann entropy for completeness, reducing attention to subsets of $\mathfrak{S}(\mathcal{H})$.\smallskip

The von Neumann entropy $H(\rho)=-\mathrm{Tr}\rho\ln\rho$
is a concave lower semicontinuous unitary invariant
function on the set $\mathfrak{S}(\mathcal{H})$ taking values in $[0,+\infty]$. It obviously has
continuous restriction to the set $\mathfrak{S}_k(\mathcal{H})$ for
each $k$. If $f=H$ then the value in the squire brackets in
(\ref{scc+}) can be expressed as follows
\begin{equation}\label{identity}
H(\rho)-\sum_i \pi_i H(\rho_i)=\sum_{i}\pi_{i}H(\rho_{i}\|\rho),
\end{equation}
 where $H(\cdot\|\cdot)$ is the quantum relative entropy defined
for arbitrary states $\rho$ and $\sigma$ in
$\mathfrak{S}(\mathcal{H})$ by the formula
$$
H(\rho\,\|\,\sigma)=\left\{\begin{array}{cc}
  \sum_{i=1}^{+\infty}\langle \varphi_i|\,\rho\ln\rho-\rho\ln\sigma\,|\varphi_i\rangle, & \mathrm{supp}\,\rho\subseteq\mathrm{supp}\,\sigma\\
  +\infty, & \mathrm{supp}\,\rho\nsubseteq \mathrm{supp}\,\sigma \\
\end{array}\right.,
$$
in which $\{|\varphi_i\rangle\}_{i=1}^{+\infty}$ is the orthonormal basis
of eigenvectors  of the operator $\rho$ (or $\sigma$) and
$\mathrm{supp}\,\rho=\mathcal{H}\ominus\mathrm{ker}\,\rho$ \cite{L,W}.
Thus we obtain from Proposition \ref{main-p} the following
continuity condition for the von Neumann entropy.\smallskip

\begin{corollary}\label{main-c-1}
\emph{The function $\rho\mapsto H(\rho)$ is continuous on a compact
subset  $\,\mathfrak{S}\subset\mathfrak{S}(\mathcal{H})$ if and
only if
\begin{equation}\label{H-cont-cond}
\lim_{k\rightarrow+\infty}\sup_{\rho\in\mathfrak{S}}\Delta^{\sigma}_k(\rho|H)=0,
\quad where \quad \Delta^{\sigma}_k(\rho|H)=\inf_{\{\pi_i,
\rho_i\}\in
M^a_\rho(\mathfrak{S}_{k}(\mathcal{H}))}\sum_{i}\pi_{i}H(\rho_{i}\|\rho).
\end{equation}}
\end{corollary}

In \cite{Sh-11} the property of an arbitrary subset
$\mathfrak{S}\subseteq\mathfrak{S}(\mathcal{H})$ expressed by
(\ref{H-cont-cond}) is called the uniform approximation property
(briefly, the UA\nobreakdash-\hspace{0pt}property) and is studied in
detail (in the extended context of the positive cone
$\mathfrak{T}_{+}(\mathcal{H})$ instead of the set
$\mathfrak{S}(\mathcal{H})$). By Corollary \ref{main-c-1} the
UA\nobreakdash-\hspace{0pt}property of an arbitrary subset $\mathfrak{S}$ is a
sufficient condition of continuity of the von Neumann entropy on
this subset and this condition is necessary if the subset $\mathfrak{S}$ is
compact.

Usefulness of the UA\nobreakdash-\hspace{0pt}property as a
continuity condition is based on possibility to analyze it by applying
well studied properties of the quantum relative entropy. This makes it possible to find a class of different
set-operations preserving the UA\nobreakdash-\hspace{0pt}property
(\cite[Proposition 4]{Sh-11}). For example,
\begin{itemize}
    \item by joint convexity and lower semicontinuity of the quantum relative entropy  the
UA\nobreakdash-\hspace{0pt}property of convex subsets
$\mathfrak{S}_1$ and $\mathfrak{S}_2$ of $\mathfrak{S}(\mathcal{H})$ implies the
UA\nobreakdash-\hspace{0pt}property of their convex closure
$\overline{\mathrm{co}}(\mathfrak{S}_1\cup\mathfrak{S}_2)$;
    \item by monotonicity of the quantum relative entropy the
UA\nobreakdash-\hspace{0pt}property a subset $\mathfrak{S}$ of $\mathfrak{S}(\mathcal{H})$ implies
the UA\nobreakdash-\hspace{0pt}property of the set
$\left\{\displaystyle \Phi(\rho)\,\left|\,\Phi\in\mathfrak{F}_n,
\rho\in\mathfrak{S}\right.\right\}$, where $\mathfrak{F}_n$ is the set
of all quantum channels having the Kraus representation consisting
of $\leq n$ summands.\footnote{The notions of a quantum channel and of
its Kraus representation are described in the next subsection.}
\end{itemize}
By using the first above-stated observation it is easy to show that
continuity of the von Neumann entropy on convex closed subsets
$\mathfrak{S}_1$ and $\mathfrak{S}_2$ of $\mathfrak{S}(\mathcal{H})$ implies its continuity on their
convex closure $\overline{\mathrm{co}}(\mathfrak{S}_1\cup\mathfrak{S}_2)$
(Corollary 7 in \cite{Sh-11}), while the second one implies the
result concerning continuity of the von Neumann entropy of
posteriory states in quantum measurements (Example 3 in
\cite{Sh-11}).

The continuity condition based on the
UA\nobreakdash-\hspace{0pt}property gives the universal method of
proving continuity of the von Neumann entropy. Various applications of this method are considered in
\cite[Section 5.2]{Sh-11}. The \emph{"only if "} part of
Corollary \ref{main-c-1} makes it possible to prove that continuity
of the von Neumann entropy on some set of states implies continuity of other important
entropic characteristics on this set (see the proofs of Corollaries
\ref{main-c-2} and \ref{main-c-3} below).

\subsection{The quantum mutual information}

Let $\mathcal{H}$ and $\mathcal{H}'$ be two separable Hilbert
spaces. A completely positive trace-preserving linear map
$\Phi:\mathfrak{T}(\mathcal{H})\rightarrow\mathfrak{T}(\mathcal{H}')$
is called \emph{quantum channel} \cite{H-SSQT,N&Ch}. By the Stinespring
dilation theorem there exist a separable Hilbert space $\mathcal{H}''$ and an
isometry $V:\mathcal{H}\rightarrow\mathcal{H}'\otimes\mathcal{H}''$
such that
\begin{equation}\label{Stinespring-rep}
\Phi(A)=\mathrm{Tr}_{\mathcal{H}''}VA V^{*},\quad \forall
A\in\mathfrak{T}(\mathcal{H}).
\end{equation}
The quantum channel
\begin{equation}\label{c-channel}
\mathfrak{T}(\mathcal{H})\ni
A\mapsto\widetilde{\Phi}(A)=\mathrm{Tr}_{\mathcal{H}'}VAV^{*}\in\mathfrak{T}(\mathcal{H}'')
\end{equation}
is called \emph{complementary} to the channel $\Phi$, it is
uniquely defined up to unitary equivalence \cite{H-c-c}. By using
representation (\ref{Stinespring-rep}) it is easy to obtain the
Kraus representation
\begin{equation}\label{Kraus-rep}
\Phi(A)=\sum_{j=1}^{+\infty}V_{j}AV^{*}_{j},\quad \forall
A\in\mathfrak{T}(\mathcal{H}),
\end{equation}
where $\{V_{j}\}_{j=1}^{+\infty}$ is a set of bounded linear
operators from $\mathcal{H}$ to $\mathcal{H}'$ such that
$\,\sum_{j=1}^{+\infty}V^{*}_{j}V_{j}=I_{\mathcal{H}}$. Via the set
$\{V_{j}\}_{j=1}^{+\infty}$ of Kraus operators of the channel
$\Phi$ its complementary channel can be expressed as follows
\begin{equation}\label{Kraus-rep-c}
\widetilde{\Phi}(A)=\sum_{i,j=1}^{+\infty}\mathrm{Tr}\left[V_{i}AV_{j}^{*}\right]|\varphi_i\rangle\langle
\varphi_j|,\quad A\in\,\mathfrak{T}(\mathcal{H}),
\end{equation}
where $\{|\varphi_i\rangle\}_{i=1}^{+\infty}$ is a particular orthonormal basis in the space
$\mathcal{H}''$ \cite{H-c-c}.

In finite dimensions ($\dim\mathcal{H},\dim\mathcal{H}'<+\infty$)
the quantum mutual information of the channel $\Phi$ at a state
$\rho\in\mathfrak{S}(\mathcal{H})$ is defined as follows
(cf.\cite{Adami})
\begin{equation}\label{mi}
 I(\rho, \Phi)=H(\rho)+H(\Phi(\rho))-H(\widetilde{\Phi}(\rho)).
\end{equation}
This is an important characteristic of a quantum channel related to
the entanglement-assisted classical capacity of this channel \cite{B&Co}.

In infinite dimensions the above definition may contain the
uncertainty $"\infty-\infty"$, but it can be modified to avoid this
problem as follows
\begin{equation}\label{mi+}
 I(\rho, \Phi) = H(\Phi \otimes \mathrm{Id}_{\mathcal{K}}
(|\varphi_{\rho}\rangle\langle\varphi_{\rho}|) \| \Phi (\rho)
\otimes \rho),
\end{equation}
where $\mathcal{K}\cong\mathcal{H}$,  $\mathrm{Id}_{\mathcal{K}}$
is the identity transformation of $\mathfrak{T}(\mathcal{H})$ and $
|\varphi_{\rho}\rangle$ is a unit vector in $\mathcal{H} \otimes
\mathcal{K}$ such that
$\mathrm{Tr}_{\mathcal{K}}|\varphi_{\rho}\rangle\langle\varphi_{\rho}|=\mathrm{Tr}_{\mathcal{H}}|\varphi_{\rho}\rangle\langle\varphi_{\rho}|=
\rho$. In \cite{H-Sh-3} it is shown that for an arbitrary quantum
channel $\Phi$ the nonnegative function $\rho\mapsto
I(\rho, \Phi)$ defined by (\ref{mi+}) is concave and lower semicontinuous on the set
$\mathfrak{S}(\mathcal{H})$ (Proposition 1) and that this function
is continuous on each subset of $\mathfrak{S}(\mathcal{H})$ on which
the von Neumann entropy is continuous, in particular, it is
continuous on the set $\mathfrak{S}_{k}(\mathcal{H})$ for each $k$
(Proposition 4). Hence for an
arbitrary quantum channel $\Phi$ the conditions of Proposition \ref{main-p}
are valid for the function $\rho\mapsto I(\rho, \Phi)$. \smallskip

By using identity (\ref{identity}), formula (\ref{mi}) and a simple
approximation it is possible to show that
\begin{equation}\label{delta-I}
\Delta_{k}^{\sigma}(\rho| I_\Phi)=\inf_{\{\pi_i, \rho_i\}\in
M^a_\rho(\mathfrak{S}_{k}(\mathcal{H}))}\sum_{i}\pi_{i}\left[H(\rho_{i}\|\rho)+H(\Phi(\rho_{i})\|\Phi(\rho))
-H(\widetilde{\Phi}(\rho_{i})\|\widetilde{\Phi}(\rho))\right],
\end{equation}
where $I_\Phi(\cdot)\doteq I(\cdot,\Phi)$,  for any state $\rho$ in
$\mathfrak{S}(\mathcal{H})$ with finite entropy. The expression in the right side of
(\ref{delta-I}) is well defined, since $\sum_{i}\pi_{i}H(\widetilde{\Phi}(\rho_{i})\|\widetilde{\Phi}(\rho))\leq \sum_{i}\pi_{i}H(\rho_{i}\|\rho)\leq H(\rho)$ \medskip
by monotonicity of the quantum relative entropy and identity (\ref{identity}).

Proposition \ref{main-p} and Corollary \ref{main-c-1} imply the
following continuity condition for the quantum mutual information,
strengthening Proposition 4 in \cite{H-Sh-3}.\medskip

\begin{corollary}\label{main-c-2}
\emph{Let $\,\Phi$ be an arbitrary quantum channel and $\,\mathfrak{S}$
be a compact subset of $\,\mathfrak{S}(\mathcal{H})$ on which the von
Neumann entropy is finite. The following assertions}
\begin{enumerate}[(i)]
    \item \emph{the function $\,\rho\mapsto H(\rho)$ is continuous on
the set $\,\mathfrak{S}$,}
\item
\emph{$\displaystyle\lim_{k\rightarrow+\infty}\sup_{\rho\in\mathfrak{S}}\Delta_{k}^{\sigma}(\rho|
I_\Phi)=0$, where $\Delta_{k}^{\sigma}(\rho| I_\Phi)$ is defined by
(\ref{delta-I}),}
\item \emph{the function $\,\rho\mapsto I(\rho, \Phi)$ is continuous on the set $\,\mathfrak{S}$,}
\end{enumerate}
\emph{are related by the implications
$\,\mathrm{(i)\Rightarrow(ii)\Leftrightarrow(iii)}$.}\smallskip

\emph{If $\,\Phi$ is a degradable channel, that is
$\widetilde{\Phi}=\Lambda\circ\Phi$ for some quantum channel
$\Lambda$, then assertions $\,\mathrm{(i)-(iii)}$ are
equivalent for an arbitrary
compact subset $\,\mathfrak{S}$ of $\,\mathfrak{S}(\mathcal{H})$.}\medskip
\end{corollary}

\textbf{Proof.} $\mathrm{(i)\Rightarrow(ii)}$ is proved by using Corollary
\ref{main-c-1}, since by monotonicity and nonnegativity of the
quantum relative entropy the expression in the square brackets in
(\ref{delta-I}) does not exceed $2H(\rho_{i}\|\rho)$
and hence $\mathrm{(ii)}$ follows from (\ref{H-cont-cond}).
$\mathrm{(ii)\Leftrightarrow(iii)}$  follows
from Proposition \ref{main-p}.\smallskip

If $\Phi$ is  a degradable channel then by using Theorem 1 in
\cite{H-Sh-3} and the 1-th chain rule from Proposition 1 in
\cite{H-Sh-3} it is easy to show that
$I(\rho,\Phi)<+\infty\Rightarrow H(\rho)<+\infty$ , while by
monotonicity of the quantum relative entropy the expression in the square
brackets in (\ref{delta-I}) is not less then
$H(\rho_{i}\|\rho)$. Thus (\ref{H-cont-cond}) follows
from $\mathrm{(ii)}$ in this case. $\square$\vspace{5pt}

\subsection{The output von Neumann entropy and the $\chi$-function of a quantum channel}

Let
$\Phi:\mathfrak{T}(\mathcal{H})\rightarrow\mathfrak{T}(\mathcal{H}')$
be a quantum channel (see Section 6.2). The output von Neumann
entropy $H(\Phi(\cdot))$ is an important characteristic involved, in
particular, in expressions for different capacities of this channel
(directly or via other characteristics) \cite{H-SSQT,N&Ch}.

The function $\rho\mapsto H_{\Phi}(\rho)\doteq H(\Phi(\rho))$ is
concave lower semicontinuous and nonnegative on the set
$\mathfrak{S}(\mathcal{H})$, but in general this function is not
continuous on sets of the family
$\{\mathfrak{S}_{k}(\mathcal{H})\}$. To apply Proposition
\ref{main-p} to the function $\rho\mapsto H_{\Phi}(\rho)$ we need
the following lemma.\medskip

\begin{lemma}\label{l-lemma}
\textit{If the function $\rho\mapsto H_{\Phi}(\rho)$ is continuous
and bounded on the set $\,\mathrm{extr}\mathfrak{S}(\mathcal{H})$
then this function is continuous on the set $\,\mathfrak{S}_{k}(\mathcal{H})$ defined by (\ref{s-k}) for each natural $\,k$.}
\end{lemma}\medskip

\textbf{Proof.} Suppose there exists a sequence
$\{\rho_{n}\}\subset\mathfrak{S}_{k}(\mathcal{H})$ converging to a
state $\rho_{0}\in\mathfrak{S}_{k}(\mathcal{H})$ such that
\begin{equation}\label{d-a}
   \lim_{n\rightarrow+\infty}H_{\Phi}(\rho_n)>H_{\Phi}(\rho_0).
\end{equation}
For each $n$ we have
$\rho_{n}=\sum_{i=1}^{k}\lambda_{i}^{n}\sigma_{i}^{n}$, where
$\{\sigma_{i}^{n}\}_{i=1}^{k}\subset\mathrm{extr}\mathfrak{S}(\mathcal{H})$
and $\{\lambda_{i}^{n}\}_{i=1}^{k}\in\mathfrak{P}_{k}$. By Lemma
\ref{g-mu-comp} we may consider that there exists
$\lim_{n\rightarrow+\infty}\lambda_{i}^{n}\sigma_{i}^{n}=A_i$ for
each $i=\overline{1,k}$, where $\{A_i\}_{i=1}^{k}$ is a set positive
trace class operators of rank $\leq1$ such that
$\rho_{0}=\sum_{i=1}^{k}A_i$. Continuity and boundedness of the
function $H_{\Phi}$ on the set $\mathrm{extr}\mathfrak{S}(\mathcal{H})$ imply
continuity of its natural extension to the cone of positive trace
class operators of rank $\leq1$ defined as follows
$$
H_{\Phi}(A) =\mathrm{Tr}A H_{\Phi}\left(\frac{A}{\mathrm{Tr}A}\right)=\mathrm{Tr}
\eta(\Phi(A)) - \eta(\mathrm{Tr}A),\quad
A\in\mathfrak{T}_{+}(\mathcal{H}),\quad \eta(x)=-x\ln x.
$$
Hence
$$
\lim_{n\rightarrow+\infty}H_{\Phi}(\lambda_{i}^{n}\sigma_{i}^{n})=H_{\Phi}(A_i),\quad
i=\overline{1,k}.
$$
By using the property of the von Neumann entropy presented after
Corollary 4 in \cite{Sh-11} we obtain a contradiction to
(\ref{d-a}). $\square$ \vspace{5pt}

The $\chi$-function of a quantum channel
$\Phi:\mathfrak{T}(\mathcal{H})\rightarrow\mathfrak{T}(\mathcal{H}')$
is a characteristic related to the classical capacity of this
channel \cite{H-SSQT,N&Ch}. It is defined as follows
$$
\chi_{\Phi}(\rho)=\sup_{\{\pi_i,\rho_i\}\in M^a(\mathfrak{S}(\mathcal{H}))}\sum_{i}\pi_{i}H(\Phi
(\rho_{i})\|\Phi(\rho)), \quad \rho\in\mathfrak{S}(\mathcal{H}).
$$
For a given subset $\mathfrak{S}$ of $\mathfrak{S}(\mathcal{H})$ the
value $\,\sup_{\rho\in\mathfrak{S}}\chi_{\Phi}(\rho)\,$ coincides with
the Holevo capacity of the $\mathfrak{S}$-constrained channel $\Phi$
\cite{Sh-3}.

The function $\rho\mapsto\chi_{\Phi}(\rho)$ is obviously concave and
nonnegative on the set $\mathfrak{S}(\mathcal{H})$. In \cite{Sh-3}
it is proved that this function is lower semicontinuous (Proposition
4) and has continuous restriction to any subset of
$\mathfrak{S}(\mathcal{H})$ on which the function $\rho\mapsto
H_{\Phi}(\rho)$ is continuous (Theorem 1). Hence Lemma \ref{l-lemma}
shows that the function $\rho\mapsto\chi_{\Phi}(\rho)$ has
continuous restriction to the set $\mathfrak{S}_{k}(\mathcal{H})$
for each $k$ if the function $\rho\mapsto H_{\Phi}(\rho)$
is continuous and bounded on the set
$\mathrm{extr}\mathfrak{S}(\mathcal{H})$.\smallskip

Thus Proposition \ref{main-p} with Lemma \ref{l-lemma} and Corollary
\ref{main-c-1} imply the following observation.\smallskip

\begin{corollary}\label{main-c-3}
\emph{Let $\,\Phi$ be a quantum channel such that the function
$\,\rho\mapsto H_{\Phi}(\rho)$ is continuous and bounded on the set
$\,\mathrm{extr}\mathfrak{S}(\mathcal{H})$. Let $\,\mathfrak{S}$ be a
compact subset of $\,\mathfrak{S}(\mathcal{H})$.  The following
assertions}
\begin{enumerate}[(i)]
    \item \emph{the function $\,\rho\mapsto H(\rho)$ is continuous on
the set $\,\mathfrak{S}$,}
\item
\emph{$\displaystyle\lim_{k\rightarrow+\infty}\sup_{\rho\in\mathfrak{S}}\Delta_{k}^{\sigma}(\rho|H_{\Phi})=0\,$,
where
$\,\displaystyle\Delta_{k}^{\sigma}(\rho|H_{\Phi})=\inf_{\{\pi_i,
\rho_i\}\in
M^a_\rho(\mathfrak{S}_{k}(\mathcal{H}))}\sum_{i}\pi_{i}H(\Phi(\rho_{i})\|\Phi(\rho))$,}
\item \emph{the function $\,\rho\mapsto H_{\Phi}(\rho)$ is continuous on
the set $\,\mathfrak{S}$,}
\item \emph{the function $\,\rho\mapsto \chi_{\Phi}(\rho)$ is continuous on
the set $\,\mathfrak{S}$,}
\end{enumerate}
\emph{are related by the implications
$\,\mathrm{(i)\Rightarrow(ii)\Leftrightarrow(iii)\Leftrightarrow(iv)}$.}\smallskip

\emph{If the output entropy of the complementary channel
$\,\widetilde{\Phi}$ is continuous on the set $\,\mathfrak{S}$ then
assertions $\,\mathrm{(i)-(iv)}$ are equivalent.}
\end{corollary}\medskip

By Corollary \ref{main-c-3} the above assertions $\mathrm{(i)-(iv)}$ are equivalent for arbitrary quantum channel $\Phi$ having Kraus representation (\ref{Kraus-rep}) with finite nonzero summands.\medskip

\textbf{Proof.} $\mathrm{(i)\Rightarrow(ii)}$ follows from Corollary
\ref{main-c-1} since monotonicity of the quantum relative entropy implies
$$
\inf_{\{\pi_i, \rho_i\}\in
M^a_\rho(\mathfrak{S}_{k}(\mathcal{H}))}\sum_{i}\pi_{i}H(\Phi(\rho_{i})\|\Phi(\rho))\leq\inf_{\{\pi_i,
\rho_i\}\in
M^a_\rho(\mathfrak{S}_{k}(\mathcal{H}))}\sum_{i}\pi_{i}H(\rho_{i}\|\rho),\quad
\forall k.
$$

$\mathrm{(ii)\Leftrightarrow(iii)}$ can be shown by applying
Proposition \ref{main-p} to the function $\rho\mapsto
H_{\Phi}(\rho)$ and by using the identity
$$
H_{\Phi}(\rho)-\sum_i \pi_i
H_{\Phi}(\rho_i)=\sum_{i}\pi_{i}H(\Phi(\rho_{i})\|\Phi(\rho)).
$$

$\mathrm{(iii)\Rightarrow(iv)}$ follows from Theorem 1 in
\cite{Sh-3}.

$\mathrm{(iv)\Rightarrow(ii)}$ can be shown by applying
Proposition \ref{main-p} to the function $\rho\mapsto
\chi_{\Phi}(\rho)$ and by using the inequality
$$
\chi_{\Phi}(\rho)-\sum_i \pi_i
\chi_{\Phi}(\rho_i)\geq\sum_{i}\pi_{i}H(\Phi(\rho_{i})\|\Phi(\rho)),
$$
valid for any $\{\pi_{i},\rho_{i}\}\in
M^a(\mathfrak{S}(\mathcal{H}))$ \cite[Proposition 4]{Sh-3}.

To prove the last assertion of the corollary it suffices to show
that continuity of the both functions $\rho\mapsto H_{\Phi}(\rho)$
and $\rho\mapsto H_{\widetilde{\Phi}}(\rho)$ on the set $\mathfrak{S}$
implies continuity of the function $\rho\mapsto H(\rho)$ on
this set.  This can be done by using Lemma \ref{old} below and
representations (\ref{Stinespring-rep}) and (\ref{c-channel}).
$\square$\medskip

\begin{lemma}\label{old}
\textit{Let $\,\{\omega_n\}$ be a sequence of states in
$\,\mathfrak{S}(\mathcal{H}\otimes\mathcal{K})$ converging to a state
$\,\omega_0$. If
$\lim_{n\rightarrow+\infty}H(\mathrm{Tr}_{\mathcal{K}}\omega_n)=H(\mathrm{Tr}_{\mathcal{K}}\omega_0)<+\infty$
and
$\,\lim_{n\rightarrow+\infty}H(\mathrm{Tr}_{\mathcal{H}}\omega_n)=H(\mathrm{Tr}_{\mathcal{H}}\omega_0)<+\infty$
then $\lim_{n\rightarrow+\infty}H(\omega_n)=H(\omega_0)<+\infty$.}
\end{lemma}\medskip

\textbf{Proof.} Let
$\omega_{n}^{\mathcal{H}}\doteq\mathrm{Tr}_{\mathcal{K}}\omega_n$
and
$\omega_{n}^{\mathcal{K}}\doteq\mathrm{Tr}_{\mathcal{H}}\omega_n$
for $n=0,1,2,...$ Since
$$
H(\omega_{n})=H(\omega_{n}^{\mathcal{H}})+H(\omega_{n}^{\mathcal{K}})-
H(\omega_{n}\|\,\omega_{n}^{\mathcal{H}}\otimes\omega_{n}^{\mathcal{K}}),
$$
we may consider that $H(\omega_{n})<+\infty$ for $n=0,1,2,...$ and
by lower semicontinuity of the quantum relative entropy we have
$$
\begin{array}{c}
\limsup\limits_{n\rightarrow+\infty}H(\omega_{n})=\lim\limits_{n\rightarrow+\infty}H(\omega_{n}^{\mathcal{H}})
+\lim\limits_{n\rightarrow+\infty}H(\omega_{n}^{\mathcal{K}})-
\liminf\limits_{n\rightarrow+\infty}H(\omega_{n}\|\,\omega_{n}^{\mathcal{H}}\otimes\omega_{n}^{\mathcal{K}})\\\\\leq
H(\omega_{0}^{\mathcal{H}})+H(\omega_{0}^{\mathcal{K}})-
H(\omega_{0}\|\,\omega_{0}^{\mathcal{H}}\otimes\omega_{0}^{\mathcal{K}})=H(\omega_{0}).
\end{array}
$$
This and lower semicontinuity of the von Neumann entropy imply
$\lim\limits_{n\rightarrow+\infty}H(\omega_{n})=H(\omega_{0})$.
$\square$

Corollary \ref{main-c-2} and Lemma \ref{a-k-closed}B imply the
following observation.\medskip
\begin{corollary}\label{main-c-3+}
\emph{Let $\,\Phi$ be a quantum channel. The following assertions
are equivalent:}
\begin{enumerate}[(i)]
    \item \emph{the function $\,\rho\mapsto H_{\Phi}(\rho)$ is
continuous and bounded on the set
$\,\mathrm{extr}\mathfrak{S}(\mathcal{H})$,}
    \item \emph{the function $\,\rho\mapsto H_{\Phi}(\rho)$ is continuous
on any subset of $\,\mathfrak{S}(\mathcal{H})$ on which the von
Neumann entropy is continuous.}
\end{enumerate}
\medskip
\end{corollary}

If assertion $\mathrm{(ii)}$ in Corollary \ref{main-c-3+} holds for a
quantum channel one can say roughly speaking that this channel
\emph{preserves continuity} of the von Neumann entropy. Assertion $\mathrm{(i)}$ in Corollary \ref{main-c-3+} can be considered as a criterion of this property.  It
implies, in particular, that the class of quantum channels
preserving continuity of the von Neumann entropy contains all quantum channels having  Kraus
representation (\ref{Kraus-rep}) with finite nonzero summands. The above criterion also shows that this class contains a quantum channel $\Phi$ if and only if it contains the complementary channel $\widetilde{\Phi}$ (since $H_{\Phi}(\rho)=H_{\widetilde{\Phi}}(\rho)$ for any $\rho\in\mathrm{extr}\mathfrak{S}(\mathcal{H})$ \cite{H-c-c}).

\section{Appendix}

\textbf{The proof of Proposition \ref{p-1}.} The function
$\hat{f}^{\mu}_{\mathcal{B}}$ is well defined on the set
$\,\overline{\mathrm{co}}(\mathcal{B})$ by Lemma \ref{bar-dec}.
Concavity of the function $\hat{f}^{\mu}_{\mathcal{B}}$ follows from
its definition and convexity of the set $M(\mathcal{B})$. By upper semicontinuity of functional
$\mathbf{f}$ defined by (\ref{functional}) and compactness of the set $M_{x}(\mathcal{B})$ for each $x$ in
$\overline{\mathrm{co}}(\mathcal{B})$
(provided by $\mu$\nobreakdash-\hspace{0pt}compactness of the set $\mathcal{A}$) the
supremum  in the definition of the value
$\hat{f}^{\mu}_{\mathcal{B}}(x)$  is achieved at a particular
measure $\mu_x$ in $M_{x}(\mathcal{B})$, that is $\hat{f}^{\mu}_{\mathcal{B}}(x)=\mathbf{f}(\mu_x)$.

Suppose the function $\hat{f}^{\mu}_{\mathcal{B}}$ is not upper
semicontinuous. Then there exists a sequence
$\{x_{n}\}\subset\overline{\mathrm{co}}(\mathcal{B})$ converging to
a point $x_{0}\in\overline{\mathrm{co}}(\mathcal{B})$ such that
\begin{equation}\label{l-s-b}
\exists\lim\limits_{n\rightarrow+\infty}\hat{f}^{\mu}_{\mathcal{B}}(x_{n})>\hat{f}^{\mu}_{\mathcal{B}}(x_{0}).
\end{equation}
As proved before for each $n$ there exists a measure $\mu_{n}\in M_{x_{n}}(\mathcal{B})$ such that
$\hat{f}^{\mu}_{\mathcal{B}}(x_{n})=\mathbf{f}(\mu_{n})$.
The $\mu$\nobreakdash-\hspace{0pt}compactness of the set $\mathcal{A}$ implies
existence of a subsequence $\{\mu_{n_{k}}\}$ converging to a
particular measure $\mu_{0}$ in $M(\mathcal{B})$. By continuity of the map
$\mu\mapsto\mathbf{b}(\mu)$ the measure $\mu_{0}$ belongs to the set
$M_{x_{0}}(\mathcal{B})$. Upper semicontinuity of
the functional $\mathbf{f}$ implies
$$
\hat{f}^{\mu}_{\mathcal{B}}(x_{0})\geq
\mathbf{f}(\mu_{0})\geq\limsup_{k\rightarrow+\infty}\mathbf{f}(\mu_{n_{k}})=
\lim_{k\rightarrow+\infty}\hat{f}^{\mu}_{\mathcal{B}}(x_{n_{k}}),
$$
contradicting to (\ref{l-s-b}).

Upper semicontinuity of the concave function
$\hat{f}^{\mu}_{\mathcal{B}}$ implies its $\mu$\nobreakdash-\hspace{0pt}concavity by Lemma
\ref{yensen}. \vspace{5pt}

\textbf{The proof of Proposition \ref{p-2}.} A) Suppose the function $\hat{f}^{\mu}_{\mathcal{B}}$ is not lower
semicontinuous. Then there exists a sequence
$\{x_{n}\}\subset\overline{\mathrm{co}}(\mathcal{B})$ converging to
a point $x_{0}\in\overline{\mathrm{co}}(\mathcal{B})$ such that
\begin{equation}\label{l-s-b+}
\exists\lim\limits_{n\rightarrow+\infty}\hat{f}^{\mu}_{\mathcal{B}}(x_{n})<\hat{f}^{\mu}_{\mathcal{B}}(x_{0}).
\end{equation}
For $\varepsilon>0$ let $\mu_{0}^{\varepsilon}$ be a measure in
$M_{x_{0}}(\mathcal{B})$ such that
$\hat{f}^{\mu}_{\mathcal{B}}(x_{0})\leq\mathbf{f}(\mu_{0}^{\varepsilon})+\varepsilon$ ($\mathbf{f}$ is the functional defined by (\ref{functional})).
By openness of the map
$M(\mathcal{B})\ni\mu\mapsto\mathbf{b}(\mu)\in\mathcal{A}$ there exists
a subsequence $\{x_{n_{k}}\}$ and a sequence
$\{\mu_{k}\}\subset M(\mathcal{B})$ converging to the
measure $\mu_{0}^{\varepsilon}$ such that
$\mathbf{b}(\mu_{k})=x_{n_{k}}$ for each $k$. Lower semicontinuity
of the functional $\mathbf{f}$ implies
$$
\hat{f}^{\mu}_{\mathcal{B}}(x_{0})\leq
\mathbf{f}(\mu_{0}^{\varepsilon})+\varepsilon\leq\liminf\limits_{k\rightarrow+\infty}
\mathbf{f}(\mu_{k})+\varepsilon\leq\lim\limits_{k\rightarrow+\infty}\hat{f}^{\mu}_{\mathcal{B}}(x_{n_{k}})+\varepsilon,
$$
contradicting to (\ref{l-s-b+}) (since $\varepsilon$ is arbitrary).

Lower semicontinuity of the concave lower bounded function
$\hat{f}^{\mu}_{\mathcal{B}}$ implies its $\mu$\nobreakdash-\hspace{0pt}concavity by Lemma
\ref{yensen}.

B) The function $\hat{f}^{\sigma}_{\mathcal{B}}$ is obviously well
defined and $\sigma$\nobreakdash-\hspace{0pt}concave on the set
$\sigma\textup{-}\mathrm{co}(\mathcal{B})$. Lower semicontinuity
of this function is proved by a simple modification of the
arguments of the proof of part A.

If
$\;\sigma\textup{-}\mathrm{co}(\mathcal{B})=\overline{\mathrm{co}}(\mathcal{B})$
then lower semicontinuity of the concave lower bounded function
$\hat{f}^{\sigma}_{\mathcal{B}}$ implies its $\mu$\nobreakdash-\hspace{0pt}concavity by Lemma
\ref{yensen}. Since
$\hat{f}^{\sigma}_{\mathcal{B}}|_{\mathcal{B}}\geq f$ by the definition of $\hat{f}^{\sigma}_{\mathcal{B}}$, we have
$$
\hat{f}^{\sigma}_{\mathcal{B}}(x)\geq \int_{\mathcal{B}}
\hat{f}^{\sigma}_{\mathcal{B}}(y) \mu(dy)\geq \int_{\mathcal{B}}
f(y)\mu(dy)
$$
for any $x\in\overline{\mathrm{co}}(\mathcal{B})$ and any measure
$\mu$ in $M_x(\mathcal{B})$. This implies
$\hat{f}^{\sigma}_{\mathcal{B}}\geq\hat{f}^{\mu}_{\mathcal{B}}$ and hence $\hat{f}^{\sigma}_{\mathcal{B}}=\hat{f}^{\mu}_{\mathcal{B}}$.
$\square$ \medskip

\textbf{The proof of Proposition \ref{basic}.} We divide the proof
into two steps.

1) Prove that for an arbitrary finitely supported measure
$\mu_0=\sum_{i=1}^{m}\pi_i \delta(x_i)$, where
$\{x_i\}_{i=1}^{m}\subset\mathcal{A}_k,
\{\pi_i\}_{i=1}^{m}\in\mathfrak{P}_m, m\in\mathbb{N}$, and an
arbitrary sequence $\{x^n\}\subset\mathcal{A}$ converging to
$x^0=\sum_{i=1}^{m}\pi_i x_i$ there exist a subsequence
$\{x^{n_k}\}$ and a sequence $\{\mu_k\}\subset M^a(\mathcal{A}_{k})$
such that $\lim_k \mu_k=\mu_0$ and $\mathbf{b}(\mu_k)=x^{n_k}$ for
all $k$.

For $k=1$ the above assertion follows from openness of the map
$M^a(\mathcal{A}_{1})\ni\mu\mapsto\mathbf{b}(\mu)$. Assume this
assertion holds for some particular $k$ and deduce its validity for
$k+1$.

Let $\mu_0=\sum_{i=1}^{m}\pi_i\delta(x_i)\in
M^a(\mathcal{A}_{k+1})$, where $\pi_i>0$ for all $i$ and
$\{x_i\}_{i=1}^{m}\nsubseteq\mathcal{A}_{k}$, let $\{x^n\}$ be a
sequence converging to $x^0=\sum_{i=1}^{m}\pi_i x_i$. For each
$i=\overline{1,m}$ we have $x_i=\alpha_i y_i+(1-\alpha_i)z_i$, where
$y_i\in\mathcal{A}_{k}$, $z_i\in\mathcal{A}_{1}$ and
$\alpha_i\in[0,1]$. Hence $x^0=\eta y^0+(1-\eta) z^0$, where
$$
\eta=\sum_{i=1}^{m}\alpha_i\pi_i\in (0,1),\quad
y^0=\eta^{-1}\sum_{i=1}^{m}\alpha_i\pi_i y_i\in \mathcal{A},\quad
z^0=(1-\eta)^{-1}\sum_{i=1}^{m}(1-\alpha_i)\pi_i z_i \in
\mathcal{A}.
$$

By stability of the set $\mathcal{A}$ we may assume (by replacing
the sequence $\{x^n\}$ by some its subsequence) existence of
sequences $\{y^n\}\subset\mathcal{A}$ and
$\{z^n\}\subset\mathcal{A}$ converging respectively to $y^0$ and
$z^0$ such that $x^n=\eta y^n+(1-\eta) z^n$. By induction
we may consider (again by passing to a subsequence) that
there exist sequences $\{\nu_n\}\subset M^a(\mathcal{A}_{k})$ and
$\{\zeta_n\}\subset M^a(\mathcal{A}_{1})$ converging to the measures
$$
\nu_0\doteq\eta^{-1}\sum_{i=1}^{m}\alpha_i\pi_i\delta(y_i)\quad
\textup{and} \quad
\zeta_0\doteq(1-\eta)^{-1}\sum_{i=1}^{m}(1-\alpha_i)\pi_i
\delta(z_i)
$$
correspondingly, such that $\mathbf{b}(\nu_n)=y^n$ and
$\mathbf{b}(\zeta_n)=z^n$ for all $n$.

By definition of the weak convergence  for arbitrary $N$ and for
arbitrary sufficiently small\footnote{In what follows it is assumed
that $\varepsilon<1/4$ and $\delta$ is so small that
$\delta$-vicinities of different points of the sets
$\{y_i\}_{i=1}^{m}$ and $\{z_i\}_{i=1}^{m}$ do not intersect each
other.} $\varepsilon>0$ and $\delta>0$ there exists such $\bar{n}>N$
that
\begin{equation}\label{rep-1}
\nu_{\bar{n}}=\sum_{i=1}^{m}\nu_{\bar{n}}^i+\nu_{\bar{n}}^r\quad
\textup{and} \quad
\zeta_{\bar{n}}=\sum_{i=1}^{m}\zeta_{\bar{n}}^i+\zeta_{\bar{n}}^r,
\end{equation}
where $\nu_{\bar{n}}^i$ and $\zeta_{\bar{n}}^i$ are measures with
finite support contained respectively in $U_\delta(y_i)$ and in
$U_\delta(z_i)$ such that
\begin{equation}\label{ineq-1}
|\nu_{\bar{n}}^i(U_\delta(y_i))-\eta^{-1}\alpha_i\pi_i|<\eta^{-1}\varepsilon\pi_i
,\quad
|\zeta_{\bar{n}}^i(U_\delta(z_i))-(1-\eta)^{-1}(1-\alpha_i)\pi_i|<(1-\eta)^{-1}\varepsilon\pi_i,
\end{equation}
all atoms of the measures $\nu_{\bar{n}}^i$ and $\zeta_{\bar{n}}^i$
have rational weights, $i=\overline{1,m}$, and
\begin{equation}\label{ineq-2}
\nu_{\bar{n}}^r(\mathcal{A})<\eta^{-1}\varepsilon, \quad
\zeta_{\bar{n}}^r(\mathcal{A})<(1-\eta)^{-1}\varepsilon.
\end{equation}
Existence of representation (\ref{rep-1}) is obvious if the sets
$\{y_i\}_{i=1}^{m}$ and $\{z_i\}_{i=1}^{m}$ consist of different
elements. If these sets contain coinciding elements existence of
this representation can be shown by "splitting" atoms of the
measures $\nu^{\bar{n}}$ and $\zeta^{\bar{n}}$ as follows. Suppose,
for example, $y_1=y_2=...=y_p=y$. Then the component
$\sum_t\lambda_t\delta(y_t)$ of the measure $\nu^{\bar{n}}$ having
atoms within $U_\delta(y)$ can be "decomposed" as
$$
\sum_t\lambda_t\delta(y_t)=\sum_t
\gamma_{1}\lambda_t\delta(y_t)+...+\sum_t
\gamma_{p}\lambda_t\delta(y_t),
$$
where $\gamma_{i}=\alpha_i\pi_i/(\alpha_1\pi_1+...+\alpha_p\pi_p)$,
and the measure $\nu_{\bar{n}}^i$ is constructed by using the
measure $\gamma_{i}\sum_t \lambda_t\delta(y_t)$.

For given $i$ let
$\nu_{\bar{n}}^i=\sum_{j=1}^{n^y_i}\frac{p^y_{ij}}{q_{i}}\delta(y_{ij})$
and
$\zeta_{\bar{n}}^i=\sum_{j=1}^{n^z_i}\frac{p^z_{ij}}{q_{i}}\delta(z_{ij})$,
where $p^*_*$ and $q_*$ are natural numbers. One can find such
natural numbers $P_{i}$, $Q^{y}_{i}$ and $Q^{z}_{i}$ that
$\sum_{j=1}^{n^y_i}\frac{p^{y}_{ij}}{q_{i}}=\frac{P_{i}}{Q^{y}_{i}}$
and
$\sum_{j=1}^{n^z_i}\frac{p^{z}_{ij}}{q_{i}}=\frac{P_{i}}{Q^{z}_{i}}$.
Let $d_i^y=(q_i Q_i^y)^{-1}$ and $d_i^z=(q_i Q_i^z)^{-1}$. By using
the "decomposition"
$$
\frac{p^y_{ij}}{q_{i}}\delta(y_{ij})=\underbrace{d_i^y\delta(y_{ij})+...+d_i^y\delta(y_{ij})}_{p^y_{ij}
Q^y_i\;\,\textup{summands}}
$$
we obtain the representation $\nu_{\bar{n}}^i=\sum_{l=1}^{P_i
q_i}d_i^y\delta(\bar{y}_i^{l})$, where $\{\bar{y}_i^{l}\}_l$ is a
set of $P_i q_i$ elements (which may be coinciding) contained in
$U_\delta(y_i)$. In the similar way we  obtain the representation
$\zeta_{\bar{n}}^i=\sum_{l=1}^{P_i q_i}d_i^z\delta(\bar{z}_i^{l})$,
where $\{\bar{z}_i^{l}\}_l$ is a set of $P_i q_i$ elements contained
in $U_\delta(z_i)$.

Let
$$
\begin{array}{c}
\displaystyle\mu_{\bar{n}}=\eta\nu_{\bar{n}}+(1-\eta)\zeta_{\bar{n}}=\sum_{i=1}^{m}(\eta\nu_{\bar{n}}^i+
(1-\eta)\zeta_{\bar{n}}^i)+\eta\nu_{\bar{n}}^r+(1-\eta)\zeta_{\bar{n}}^r\\=
\displaystyle\sum_{i=1}^{m}\sum_{l=1}^{P_i q_i}\left(\eta
d_i^y\delta(\bar{y}_i^{l})+(1-\eta)d_i^z\delta(\bar{z}_i^{l})\right)+\eta\nu_{\bar{n}}^r+(1-\eta)\zeta_{\bar{n}}^r
\end{array}
$$
be a measure with the barycenter $\eta y^{\bar{n}}+(1-\eta)z^{\bar{n}}=x^{\bar{n}}$. The
measure
$$
\hat{\mu}_{\bar{n}}=\sum_{i=1}^{m}\sum_{l=1}^{P_i q_i}(\eta
d_i^y+(1-\eta)d_i^z)\delta(\bar{x}_i^l)+\eta\nu_{\bar{n}}^r+(1-\eta)\zeta_{\bar{n}}^r,
\;\,\textup{where}\;\,\bar{x}_i^l=\frac{\eta
d_i^y\bar{y}_i^{l}+(1-\eta)d_i^z\bar{z}_i^{l}}{\eta
d_i^y+(1-\eta)d_i^z}
$$
has the same barycenter and lies in $M^a(\mathcal{A}_{k+1})$. Since
$$
\bar{\alpha}_i=\frac{\eta d_i^y}{\eta
d_i^y+(1-\eta)d_i^z}=\frac{\frac{\eta P_i}{Q^y_i\pi_i}}{\frac{\eta
P_i}{Q^y_i\pi_i}+\frac{(1-\eta)P_i}{Q^z_i\pi_i}}
$$
and (\ref{ineq-1}) implies $|\frac{\eta
P_i}{Q^y_i\pi_i}-\alpha_i|<\varepsilon$,
$|\frac{(1-\eta)P_i}{Q^z_i\pi_i}-(1-\alpha_i)|<\varepsilon$, it
is easy to show that $|\bar{\alpha}_i-\alpha_i|<6\varepsilon$. Thus
we conclude that
$\bar{x}_i^l=\bar{\alpha}\bar{y}_i^l+(1-\bar{\alpha})\bar{z}_i^l\in
U_{\delta(i)}(x_i)$ for all $i=\overline{1,m}$ and $l=\overline{1,P_i
q_i}$, where $\delta(i)=2\delta+C_{y_i,z_i}(6\varepsilon)$ (see
Remark \ref{perturb}). Since $P_i q_i(\eta
d_i^y+(1-\eta)d_i^z)=\eta\frac{P_i}{Q^y_i}+(1-\eta)\frac{P_i}{Q^z_i}$,
by using (\ref{ineq-1}) and (\ref{ineq-2}) it is easy to show that
\begin{equation}\label{ineq-3}
|\hat{\mu}_{\bar{n}}(U_{\delta(i)}(x_i))-\pi_i|\leq 4\varepsilon
\end{equation}
provided $U_{\delta(i)}(x_i)\cap U_{\delta(i')}(x_{i'})=\emptyset$ for all
$i\neq i'$.

For natural $l$ let $n_l=\bar{n}$ and $\mu_l=\hat{\mu}_{\bar{n}}$,
where $\bar{n}$ and $\hat{\mu}_{\bar{n}}$ are produced by the above
construction with $N=l$ and $\varepsilon=\delta=1/l$. Then
$\mathbf{b}(\mu_l)=x^{n_l}$ and (\ref{ineq-3}) implies weak
convergence of the sequence $\{\mu_l\}$ to the measure $\mu_0$.\smallskip

2) Let $\mu_0=\sum_{i=1}^{+\infty}\pi_i \delta(x_i)$ be an arbitrary
measure in $M^a(\mathcal{A}_k)$ and $\,\{x^n\}\subset\mathcal{A}\,$
be a sequence converging to
$x^0=\sum_{i=1}^{+\infty}\pi_i x_i$. For natural $m$ let
$\mu^m_0=(\lambda_m)^{-1}\sum_{i=1}^{m}\pi_i \delta(x_i)$, where
$\lambda_m=\sum_{i=1}^{m}\pi_i$, and let $\nu^m_0$ be a
measure in $M^a(\mathcal{A}_1)$ such that
$\mathbf{b}(\nu^m_0)=(1-\lambda_m)^{-1}\sum_{i>m}\pi_i x_i$.

Since the sequence $\{\mu^m_0\}_m$ converges to the measure $\mu_0$,
for given natural $l$ there exists $m_l$ such that $\mu^{m_l}_0\in
U_{1/l}(\mu_0)$ and $\lambda_{m_l}>1-1/l$.\footnote{The set $M(\mathcal{A})$ can be considered as a metric space \cite{Par}.} We have
$x^0=\lambda_{m_l}\mathbf{b}(\mu^{m_l}_0)+(1-\lambda_{m_l})\mathbf{b}(\nu^{m_l}_0)$.
By stability of the set $\mathcal{A}$ we may assume (by replacing
the sequence $\{x^n\}$ by some its subsequence) existence of
sequences $\{y^n\}\subset\mathcal{A}$ and
$\{z^n\}\subset\mathcal{A}$ converging respectively to
$\mathbf{b}(\mu^{m_l}_0)$ and $\mathbf{b}(\nu^{m_l}_0)$ such that
$x^n=\lambda_{m_l}y^n+(1-\lambda_{m_l})z^n$.

By the first part of the proof we may consider (again by passing to
a subsequence) that there exists a sequence $\{\mu_n\}\subset
M^a(\mathcal{A}_k)$ converging to the measure $\mu^{m_l}_0$ such
that $\mathbf{b}(\mu_n)=y_{n}$ for all $n$. Hence there exists $n_l>l$
such that $\mu_{n_l}\in U_{1/l}(\mu^{m_l}_{0})\subset
U_{2/l}(\mu_{0})$. Let
$$
\bar{\mu}_l=\lambda_{m_l}\mu_{n_l}+(1-\lambda_{m_l})\nu_{n_l},
$$
where $\nu_{n_l}$ is an arbitrary measure in $M^a(\mathcal{A}_1)$
such that $\mathbf{b}(\nu_{n_l})=z^{n_l}$.

It is easy to see that the sequence $\{\bar{\mu}_l\}$ is contained
in $M^a(\mathcal{A}_{k})$ and converges to the measure $\mu_0$ while
by the construction $\mathbf{b}(\bar{\mu}_l)=\lambda_{m_l}y^{n_l}+(1-\lambda_{m_l})z^{n_l}=x^{n_l}$ for each $l$.
$\square$\bigskip

I am grateful to A.S.Holevo and the participants of his seminar for useful discussion.
I am also grateful to the organizers of the
workshop \emph{Thematic Program on Mathematics in Quantum
Information} at the Fields Institute, where some of this work was done.

\end{document}